\documentclass[11pt]{amsart}
\usepackage{ifthen, citesort}
\usepackage{amssymb}
\usepackage[centertags]{amsmath}
\usepackage{amscd}
\usepackage{amsxtra}
\usepackage{latexsym}
\newtheorem{thm}{Theorem}[section]
\newtheorem{prop}[thm]{Proposition}
\newtheorem{lem}[thm]{Lemma}
\newtheorem{cor}[thm]{Corollary}
\newtheorem{definition}[thm]{Definition}
\newtheorem{remark}[thm]{Remark}
\newtheorem{example}[thm]{Example}

\newenvironment{rem}{\begin{remark}\begin{rm}}{\end{rm}\end{remark}}
\newenvironment{defi}{\begin{definition}\begin{rm}}{\end{rm}\end{definition}}

\newcommand{\bt}{\begin{thm}}
\newcommand{\br}{\begin{rem}}
\newcommand{\bc}{\begin{cor}}
\newcommand{\bl}{\begin{lem}}
\newcommand{\bd}{\begin{defi}}
\newcommand{\bp}{\begin{prop}}

\newcommand{\et}{\end{thm}}
\newcommand{\er}{\end{rem}}
\newcommand{\ec}{\end{cor}}
\newcommand{\el}{\end{lem}}
\newcommand{\ed}{\end{defi}}          
\newcommand{\ep}{\end{prop}}

\def\qedbox{\mbox{\large $\Box$}}
\newenvironment{Pf}{\noindent{\it
       Proof:$\;\; $}}{\rule{0mm}{0mm}\nolinebreak\hfill \qedbox \vspace{12pt}}
\newcommand{\bpf}{\begin{Pf}}
\newcommand{\epf}{\end{Pf}}

\newcommand{\boxma}{\\\begin{minipage}[b]{14cm}}
\newcommand{\boxme}{\end{minipage} \hfill $\Box$}
\newcommand{\mb}[1]{\mathbb{#1}}
\newcommand{\rn}{\ensuremath{\mb{R}^{n}}}
\newcommand{\hn}{\ensuremath{\mathcal{H}^{n}}}
\newcommand{\hne}{\ensuremath{\mathcal{H}^{n+1}}}
\newcommand{\hnz}{\ensuremath{\mathcal{H}^{n+2}}}
\newcommand{\hz}{\ensuremath{\mathcal{H}^{2}}}

\newcommand{\R}{\ensuremath{\mb{R}}}
\newcommand{\N}{\mb{N}}

\newcommand{\rne}{\ensuremath{\mb{R}^{n+1}}}
\newcommand{\rnz}{\ensuremath{\mb{R}^{n+2}}}
\newcommand{\res}{\ensuremath{\,\textsf{\small \upshape L}\,}}
\newcommand{\p}{\partial}
\newcommand{\beq}{\begin{equation}{}}
\newcommand{\eeq}{\end{equation}{}}
\newcommand{\beqs}{\begin{equation*}}
\newcommand{\eeqs}{\end{equation*}}
\newcommand{\lb}[1]{\label{#1}}
\newcommand{\dt}{\ensuremath{\frac{\p}{\p t}}}
\newcommand{\eps}{\varepsilon}
\newcommand{\net}{N^\eps_t}
\newcommand{\neit}{N^{i}_t}
\newcommand{\meit}{\mu^{\eps_i}_t}
\newcommand{\inl}{\int\limits}

\renewcommand{\leq}{\leqslant}
\renewcommand{\geq}{\geqslant}
\newcommand{\ra}{\rightarrow}

\begin{document}

\title[Nonlinear Evolution by Mean Curvature]{Nonlinear Evolution by Mean Curvature and Isoperimetric Inequalities}

\author{Felix Schulze}
\address{Felix Schulze, FU Berlin, Arnimallee 6,  14195 Berlin, Germany}
\email{Felix.Schulze@math.fu-berlin.de}
\thanks{The author is a member of SFB 647/B3 ``{Raum - Zeit - Materie}''}


\begin{abstract}
\noindent 
Evolving smooth, compact hypersurfaces in $\R^{n+1}$ with normal speed
equal to a positive power $k$ of the mean curvature improves a certain
'isoperimetric difference' for $k\geq n-1$. As singularities may develop
before the volume goes to zero, we develop a weak level-set
formulation for such flows and show that the above monotonicity 
is still valid. This proves the isoperimetric inequality for $n \leq 7$. 
Extending this to complete, simply connected $3$-dimensional manifolds 
with nonpositive sectional curvature, we give a new proof for the
Euclidean isoperimetric inequality on such manifolds.
\end{abstract}

\maketitle


\section{Introduction}
 Let $M^n$ be a smooth $n$-dimensional
compact manifold without boundary and \mbox{$F_0:M^n\ra N^{n+1}$} a
smooth embedding into an $n\!+\!1$-dimensional Riemannian manifold
$(N^{n+1}, \bar{g})$. We assume further that $F_0(M)$ has positive
mean curvature in $N^{n+1}$. Starting from such an initial
hypersurface there exists, at least for a short time interval $[0,T)$,
an evolution \mbox{$F(\cdot,t):M^n\times [0,T)\ra N^{n+1}$}, which
satisfies
\begin{equation}\tag{\ensuremath{\star}}\begin{cases}
\ \ \, F(\cdot , 0) = F_0(\cdot) \\[+1.0ex]
\ \displaystyle{\frac{dF}{dt}} (\cdot, t) = -H^k (\cdot, t) \nu(\cdot ,t)\ .
\end{cases}\qquad\qquad\qquad\qquad\qquad\qquad\qquad\qquad
\end{equation}

where $k \geq 1,\ H$ is the mean curvature and $\nu$ is the outer unit normal, such
that $-H\nu=\mathbf{H}$ is the mean curvature vector. Let $A(t)$ denote
the area of such an evolving hypersurface, $V(t)$ the enclosed volume,
and $c_{n+1}$ the Euclidean
isoperimetric constant. We aim to exploit the following fact, which 
G. Huisken has drawn our attention to: the 'isoperimetric difference'
\beq\lb{intro.1} A(t)^{\frac{n+1}{n}}- c_{n+1}V(t)
\eeq
is monotonically decreasing under such a flow, provided $k\geq n-1$
and the inequality
\beq\lb{intro.2}
\int_{M_t}|\mathbf{H}|^n\, d\mu \geq \Big(\frac{n}{n+1}\, c_n\Big)^n
\eeq 
holds on all of the evolving surfaces for $t\in (0,T)$. In the case
that $N=\R^{n+1}$ an easy calculation proves this inequality for an
arbitrary closed hypersurface which is at least $C^2$. If $n=2$ and $N^3$ has negative
sectional curvatures we can use the monotonicity formula to show that
(\ref{intro.2}) holds on any closed hypersurface. 

\medskip

If we assume that the flow $(\star)$ exists until $V(t)$ decreases to
zero, the monotonicity of $(\ref{intro.1})$ would prove the Euclidean
isoperimetric inequality for this initial
configuration. Unfortunately, without special geometric assumptions
(see \cite{Huisken84},\cite{convex}),
singularities may develop even before the volume goes to zero.
To cope with this problem, we replace $(\star)$ by the following
level-set formulation. Let $\Omega \subset N$ be a bounded,
open set with smooth boundary $\p\Omega$, such that $\p\Omega$ has
positve mean curvature. The evolving surfaces are then given as
level-sets of a
continuous function $u:\bar{\Omega} \ra \R^+,\ u=0$ on $\p\Omega$ via
$$ \Gamma_t=\p \{x \in \Omega\ |\ u(x)> t\},$$
and $(\star)$ is replaced by the degenerate elliptic equation
\beq\tag{\ensuremath{\star\star}}
\text{div}_N\bigg(\frac{\nabla u}{|\nabla u|}\bigg) = -\frac{1}{|\nabla
  u|^{\frac{1}{k}}}\ .
\eeq
Here the left hand side gives the mean curvature of the level-sets and
the right-hand side is the speed, raised to the appropriate
power. If $u$ is smooth at $x\in\Omega$ with nonvanishing gradient,
then $(\star\star)$ implies that the level-sets $\Gamma_t$ are
evolving smoothly according to $(\star)$ in a neighborhood of $x$.
This formulation is inspired by the work of Evans-Spruck
\cite{EvansSpruckI} and Chen-Giga-Goto \cite{ChenGigaGoto} on mean
curvature flow and by the work of Huisken-Ilmanen
\cite{HuiskenIlmanen} on the inverse mean curvature flow.

We use the method of elliptic regularisation to define a family of
approximating problems to $(\star\star)$. We prove the existence of
smooth solutions to the approximating problem, which by a
uniform a-priori gradient bound
subconverge to a lipschitz continuous function $u$ on
$\Omega$. We define any such limit function $u$ to be a weak solution
to $(\star\star)$ and call it a weak $H^k$-flow generated by
$\Omega$. 
This is justified since such a weak solution
solves $(\star\star)$ in the viscosity sense and even more we show
that as long as the smooth solution to $(\star)$ exists, it coincides
with any weak solution. In the case that $n\leq 6$ and the ambient
space is flat we can show that this weak solution is unique, i.e. it
does not depend on the approximating sequence $u^{\eps_i}$.

\bt \lb{mo.t1}Let $\Omega$ be a bounded, open subset
of $N$, such that $H|_{\p\Omega}>0$. If $n=2$ let $k\geq 1$
and $N$ be a complete, simply connected 3-manifold with
nonpositive sectional curvatures. If $n\geq 3$ let
$N=\rne$ and $k>n$. If $u$ is a weak $H^k$-flow generated by
$\Omega$, then the isoperimetric difference 
$$ I_t:= \big(\hn(\p^*\{u>t\})\big)^\frac{n+1}{n} - c_{n+1}\hne(\{u>t\}) $$
is a nonnegative, monotonically decreasing function on $[0,T)$, where
$T:=\sup_\Omega u$.
\et

Here we denote with $\mathcal{H}^l$ the $l$-dimensional Hausdorff-measure.
We can use this monotone quantity to give a proof of the isoperimetric
inequality in $\rne$ for $n\leq 7$. The same technique also works if
the ambient manifold is simply connected and complete with nonnegative
sectional curvatures, which gives a new proof of the result by B. Kleiner in \cite{Kleiner92}.

\bc[{\bf Isoperimetric inequality}]\lb{mo.c1} Let $U \subset \rne$ be
a compact domain with smooth boundary and $n+1\leq 8$, or $U\subset
N^3$, where $N^3$ is as above. Then 
\beq\lb{mo.2}
\big(\hn(\p U)\big)^\frac{n+1}{n} \geq c_{n+1}\hne(U) \ .
\eeq
\ec
\vspace{2ex}
In the case that $N^3$ has sectional curvatures bounded
above by $-\kappa$, $\kappa>0$, define the function $f:\R^+\ra \R^+$ by
\beq\lb{intro.3} f_\kappa(A):=\int_0^A\frac{a^\frac{1}{2}}{(16\pi+ 4\kappa a)^\frac{1}{2}}\,
  da \ .\eeq
Let $\Omega\subset N^3$ be open and bounded and $u$ be weak
$H^k$-flow generated by $\Omega$. We show that under 
the restriction that all superlevelsets $\{u>t\}$ minimize area from
the outside in $N$ the quantity
$$I_t^\kappa:= f_\kappa(\hz(\p^*(\{u>t\})) - \mathcal{H}^3(\{u>t\}) $$
is a nonnegative, monotonically decreasing function for  $t\in
[0,T)$. This enables us to give a proof of the following stronger result, which also already
appeared in \cite{Kleiner92}.

\bt\lb{mo.t2}Let $N^3$ be a complete, simply connected,
3-dimensional Riemannian manifold with sectional curvatures bounded
above by $-\kappa <0$. If $U\subset N^3$ is a compact domain
with smooth boundary, then
$$f_\kappa\big(\hz(\p U)\big)\geq \mathcal{H}^3(U)\ .$$
Moreover, equality holds for geodesic balls in the model space
$N^3_\kappa$ with constant sectional curvature $-\kappa$.
\et

Mean curvature flow in the level set formulation was developed in
\cite{EvansSpruckI} and \cite{ChenGigaGoto}, see also
\cite{Ilmanen}. G. Huisken and T. Ilmanen developed a weak level set formulation
for the inverse mean curvature flow to prove the Riemannian Penrose
inequality in \cite{HuiskenIlmanen}. \\
As already mentioned before, B. Kleiner proved the Euclidean
isoperimetric inequality on a complete, simply connected 3-manifold
with nonpositive sectional curvatures in \cite{Kleiner92}. In the case
that the ambient manifold is 4-dimensional the corresponding result was proven by C. Croke,
\cite{Croke84}.\\
The monotonicity of (\ref{intro.1}) under mean curvature flow for
$n=1,2$ was also observed by P. Topping. In \cite{Topping98} he uses the monotonicity under curve
shortening flow to prove optimal
isoperimetric inequalities on 2-surfaces. Utilizing the monotonicity
 under mean curvature flow of 2-surfaces he gives sufficient geometric 
conditions for the formation of singularities
under this flow.\\
Similar monotonicities of area and volume 
under the affine normal flow were employed by B. Andrews in
\cite{Andrews96} to give new proofs of affine isoperimetric
inequalities.\\[2ex]
{\bf Outline.} In \S2, we show how to derive the monotonicity of
(\ref{intro.1}) in case that the flow is smooth.
In \S3 we define by
elliptic regularisation the $\eps$-regularised version of
$(\star\star)$. Using barrier techniques we
prove uniform a-priori sup and gradient bounds, which we apply to
show existence of solutions $u^\eps$ to the regularised problem. Also by these
a-priori bounds the solutions $u^\eps$ subconverge as $\eps \ra 0$ to
a lipschitz-continuous function $u$ on $\Omega$ which we define to be
a weak solution to $(\star\star)$. \\
In \S4 we establish that a weak solution satisfies an 
avoidance principle w.r.t.\! smooth $H^k$-flows. 
For flat ambient space with $n+1\leq 7$ we use this to prove uniqueness.\\
The approximating solutions $u^\eps$ have the important geometric
property that, scaled appropriately, they constitute a smooth, graphical,
translating solution to the $H^k$-flow in $\Omega\times \R$. This can
be applied to obtain an approximation of the weak $H^k$-flow by
smooth flows in one dimension higher.\\
In \S5 we refine our
understanding of this approximation and use it to prove properties of
the weak limit flow. These properties include that the weak flow is
non-fattening, i.e. the sets $\{u=t\}\subset \rne$ do not develop positive
$\hne$-measure. We also show that the sets $\{u>t\}$ minimize area
from the outside in $\Omega$. As another consequence the level sets
$\{u=t\}$ are actually quite nice, i.e. for $k> n-1$ and almost every
$t$ they are $C^{1,\alpha}$-hypersurfaces up to a closed set of $\hn$-measure
zero.\\
In \S6, we prove that the estimate (\ref{intro.2}) holds on $\Gamma_t$
for a.e. $t$, provided $k>n$ and the ambient space is flat. 
To do this, we first replace $\Gamma_t$ by
an outer equidistant hypersurface to the convex hull of $\{u>t\}$. Since
such a hypersurface is convex and $C^{1,1}$ we can apply the same
proof as in the smooth case to show (\ref{intro.2}) on this
hypersurface. The biggest chunk of work in this chapter is then to translate this
estimate back to $\Gamma_t$. The main ingredient
there is an estimate on the growth of the area of the equidistant
hypersurfaces in terms of an integral of the mean curvature of
$\Gamma_t$. 
Here again we use the approximation by smooth translating flows in $\Omega \times
\R$. In the case that $n=2$ and the ambient space is not flat, we give
a proof of (\ref{intro.2}) which uses the monotonicity formula and
thus needs much less regularity of $\Gamma_t$. The stronger estimate
needed for Theorem \ref{mo.t2} is proved by combining the techniques
in the flat case with the Gauss-Bonnet formula. \\
Finally in \S7 we use the approximation by smooth flows in $\Omega
\times \R$ together with the estimate from \S6 and a lower
semicontinuity argument to show that the monotonicity of (\ref{intro.1})
holds in the limit. This is then applied to yield a proof of the stated
isoperimetric inequalities. The restriction on the dimension of the
ambient space in the flat case comes from the problem that to start
the flow, we have to replace a bounded set $U\subset \rne$ with smooth
boundary by its outer minimizing hull, which is only known to be smooth,
more precisely $C^{1,1}$, for $n\leq 6$. For $n>6$ the
outer minimizing hull can have singularities on the part
away from the obstacle $U$. For $n=7$ these singularities
are still isolated and
an argument of R. Hardt and L. Simon can be applied to show that we can perturb $U$
slightly such that its outer minimizing hull again is $C^{1,1}$. \\[2ex]
We want to especially 
thank G. Huisken for many stimulating discussions and support. 
We also want to thank K. Ecker and T. Ilmanen for further
discussions and support. Finally we want to thank B. White for bringing to our attention the
argument of Hardt-Simon above.


\section{The smooth Case}\lb{smooth}

\noindent Given a smooth solution to $(\star)$, not necessarily compact, we
first compute the evolution equations of geometric quantities like 
the induced metric $g_{ij}$, the induced measure $d\mu$ and the mean curvature $H$. 

\bl \lb{evolequ}
The following evolution equations hold.

\begin{itemize}
\item[i)] $\frac{\partial}{\partial t}\, g_{ij}=-2H^k h_{ij}\ ,$
\item[ii)] $\frac{\partial}{\partial t}\, H=kH^{k-1}\Delta H
+k(k-1)H^{k-2}|\nabla H|^2+|A|^2H^k +\overline{\text{R}}\text{ic}(\nu,\nu)H^k \ ,$
\item[iii)] $\frac{\partial}{\partial t}\,d\mu = -H^{k+1}d\mu\ $
\end{itemize}
\el 
{\it Proof:}\ All of the above follows from a direct calculation as
for example  in \cite{Huisken84}.\\[+2ex]
In the case that we have a smooth solution of closed hypersurfaces in $\rne$ to
$(\star)$, we first demonstrate how to show the monotonicity of 
(\ref{intro.1}) as claimed in the introduction. Here
$$ c_{n+1}= \big((n+1)^{n+1}\omega_{n+1}\big)^\frac{1}{n} $$
is the Euclidean isoperimetric constant, $\omega_{n+1}$ denoting the
volume of the unit ball in $\rne$. To do this, let us first prove
estimate (\ref{intro.2}) for an arbitrary closed hypersurface
$M\subset \rne$ which is at least $C^{1,1}$. Let $M^+$ be the intersection of $M$ with the
boundary of its outer convex hull. Since the unit normal map $\nu$ (let us always
choose the outer unit normal), restricted to $M^+$, covers
$\mathbb{S}^n$ at least once we can estimate
\beq\lb{sm.2}
|\mathbb{S}^n| \leq \inl_{M^+}\nu^*\!do_{S^n}= \inl_{M^+}G\, d\hn \leq
\frac{1}{n^n}\inl_{M^+}H^n\, d\hn \leq \frac{1}{n^n}\inl_{M}|\mathbf{H}|^n\,
d\hn \ ,
\eeq
which is (\ref{intro.2}). Here $G$ denotes the Gauss curvature and we used
that on $M^+$ all principal curvatures are nonnegative, thus we can apply the
arithmetic-geometric mean inequality in the second estimate. Note that
for $k\geq n-1$ this implies by H\"older
\beq\lb{sm.3}
n^n (n+1)\,\omega_{n+1}\leq \bigg(\inl_M
|H|^{k+1}d\hn\bigg)^\frac{n}{k+1}|M|^{1-\frac{n}{k+1}}\ .
\eeq
Now use the evolution equations and the above estimate to calculate
\beq\lb{sm.4}\begin{split}
-\frac{d}{dt} V =&\ \inl_{M_t}H^k\, d\hn \leq \bigg(\inl_{M_t}H^{k+1}\,
d\hn\bigg)^\frac{k}{k+1}A^\frac{1}{k+1}\\
&\ \cdot\frac{1}{n}\big((n+1)\omega_{n+1}\big)^{-\frac{1}{n}}
\bigg(\inl_{M_t}H^{k+1}\,
d\hn\bigg)^\frac{1}{k+1}A^{\frac{1}{n}-\frac{1}{k+1}}\\
\leq&\ \frac{1}{n}\big((n+1)\omega_{n+1}\big)^{-\frac{1}{n}} \inl_{M_t}H^{k+1}\,
d\hn \cdot A^\frac{1}{n}= -\frac{1}{c_{n+1}}\frac{d}{dt}A^\frac{n+1}{n}\ .
\end{split}
\eeq
Rearranging, this implies $\frac{d}{dt} I(t)\leq 0\ $.\\[2ex]
If the ambient manifold $N$ is 3-dimensional and has nonpositve sectional
curvatures, it needs some more work to prove (\ref{intro.2}). If
$M\subset N$ is a closed hypersurface which is diffeomorphic
to a sphere and at least $C^{1,1}$ one can use the Gauss-Bonnet
formula, see (\ref{me.34}). In Lemma \ref{me.l6} we give a
proof which works without any restriction on the topology. The proof uses a
variant of the monotonicity formula, and thus needs much less
regularity of $M$. The calculation (\ref{sm.4}) also applies in this
setting to show that $I(t)$ is decreasing in time.


\section{Elliptic regularisation}\lb{ellipticreg}

To define a weak solution of
$(\star\star)$ we apply an approximation scheme known as elliptic
regularisation. Similar techniques to show the existence of weak
solutions, often in the viscosity sense, have been used by various
authors, see \cite{EvansSpruckI}, 
\cite{HuiskenIlmanen}, \cite{Ilmanen}. We define the following
approximating equation.
\begin{equation*}
 (\star\star)_\eps \begin{cases}
 \ \ \text{div}_{N}\Big(\frac{\nabla u^\eps}
{\sqrt{\eps^2+|\nabla u^\eps|^2}}
\Big ) =-\big(\eps^2+|\nabla u^\eps|^2\big)^{-\frac{1}{2k}}\ \ 
&\text{in} \ \Omega\\[+2.0ex]
\ \ u^\varepsilon = 0 &\text{on}\ \p\Omega
\end{cases}\qquad\qquad\qquad
\end{equation*}
We can give this equation a geometric interpretation. It implies that
the downward translating graphs
$$
\net:=\text{graph}\bigg(\frac{u^\eps(x)}{\eps}-\frac{t}{\eps}\bigg),\qquad
-\infty<t<\infty $$
solve the $H^k$-flow $(\star)$ smoothly in the manifold
$\Omega\times \R$. To verify this, define the function
\beq\lb{er.1}
U^\eps(x,z):=u^\eps(x) - \eps z,\qquad (x,z)\in \Omega\times\R, 
\eeq
such that $\{U^\eps=t\}=\net$. If we assume smoothness, one can check
that $U^\eps$ satisfies $(\star\star)$ on $\Omega\times\R$ if and only if
$u^\eps$ satisfies $(\star\star)_\eps$ on $\Omega$.\\[2ex]
Let us now assume that the solutions $u^\eps$ converge in a suitable
sense to a weak solution $u$ with level sets $\{u=t\}$. The geometric
idea in this approximation then is that the possibly singular
evolution of $\{u=t\}$ is well approximated by the evolution of $\net$
in the sense that $\net \approx \{u=t\}\times \R$ for sufficiently
small $\eps >0$.\\
To show the existence of solutions to
$(\star\star)_\eps$ we first have to prove a-priori sup- and
gradient-bounds.

\bl\lb{er.l1} Let $u^\eps$ be a smooth solution to
$(\star\star)_\eps$. Then for  $0<\eps\leq 1$, 
\beq\lb{er.2}
\sup_\Omega |u^\eps| \leq C(n,k,\text{diam}(\Omega))\ .
\eeq
\el
\bpf Our aim is to construct a supersolution. So pick a $p_0 \in
N$ with dist$(p_0,\Omega)=1$. Let $S_r:= \p B(p_0,r)$. Since
$N$ is diffeomorphic via the exponential map at $p_0$ to
$\rne$, the hypersurfaces $S_r$ are smooth. The evolution of the mixed
second fundamental form is given by
$$\frac{\p}{\p r}h^i_{\ j} = -h^i_{\ k}h^k_{\ j} - R_{0\ 0j}^{\ i}\geq
-h^i_{\ k}h^k_{\ j},
$$
since we have assumed nonpositve sectional curvatures. Since $S_r$ is
convex for small $r$, this implies that this remains so for all $r>0$. Taking the
trace, we see that the mean curvature $H$ of $S_r$ satisfies
$$\frac{\p}{\p r} H = -|A|^2-\text{Ric}(\nu,\nu) \geq -H^2.$$
Using that $\lim_{r\ra 0} H = +\infty$ and integration  implies
\beq\lb{er.3}
H(p)\geq \frac{1}{r}\ .
\eeq 
for $p\in S_r$. We make the ansatz $\Phi(p)=\psi(r)$, where
$r(x)=\text{dist}(x,p_0)$, and compute
\beqs
\begin{split}
\text{div}_N\bigg(\frac{\nabla \Phi}{\sqrt{\eps^2+|\nabla
    \Phi|^2}}\bigg) =
&\ \frac{\psi^\prime}{\sqrt{\eps^2+(\psi^\prime)^2}}\, \Delta r +
\bar{g}\big(\nabla\bigg(\frac{\psi^\prime}{\sqrt{\eps^2+(\psi^\prime)^2}}\bigg),\nabla
r\big)\\
=&\  \frac{\psi^\prime}{\sqrt{\eps^2+(\psi^\prime)^2}}\, H(p)+
\frac{\eps^2\psi^{\prime\prime}}{\big(\eps^2+(\psi^\prime)^2\big)^{3/2}}\ .
\end{split}
\eeqs 
which should be less than $-(\eps^2 +(\psi^\prime)^2)^{-1/(2k)}$ for
$\Phi$ to be a supersolution. Let us assume that $\psi^\prime\leq
0$. We apply (\ref{er.3}) to see that a sufficient condition is that
\beq\lb{er.5}
\frac{1}{r}\geq -
\frac{\eps^2\psi^{\prime\prime}}{\psi^\prime\big(\eps^2+(\psi^\prime)^2\big)}
-\frac{1}{\psi^\prime}\big(\eps^2+(\psi^\prime)^2\big)^{\frac{k-1}{2k}}\
.
\eeq

Now assume $\Omega \subset B(p_0,R_0)$ for some $R_0$ large
enough. Let $\sigma>0$ be a constant still to be chosen and take $\psi=
\sigma (k+1)^{-1} (R_0^{k+1}-r^{k+1})$, which gives 
$$ \psi^\prime = -\sigma r^k, \qquad \psi^{\prime\prime}= -\sigma k
r^{k-1}\ .$$
The inequality (\ref{er.5}) then becomes 
$$
\frac{1}{r}\geq - \frac{\eps^2 k}{r(\eps^2 +\sigma ^2
  r^{2k})}+\frac{1}{\sigma r^k} (\eps^2+\sigma^2r^{2k})^\frac{k-1}{2k}
\ .
$$
Dropping the first term on the RHS, a sufficient condition again is
$$ 1 \geq \frac{1}{\sigma}\Big(\frac{\eps^2}{r^2k}
+\sigma^2\Big)^\frac{k-1}{2k} = \frac{1}{\sigma^\frac{1}{k}}\Big(\frac{\eps^2}{\sigma^2r^{2k}}
+1\Big)^\frac{k-1}{2k}\ .
$$
Since $r\geq 1$ on $\Omega$, we can choose $\sigma$ large enough such
that the above condition is satisfied for all $p\in \Omega$ and
$\eps\leq 1$. Thus $\Phi$ is a positive supersolution on $\Omega$ for all
$0<\eps\leq 1$. By the maximum principle we obtain the desired
sup-bound. Since every solution to $(\star\star)_\eps$ has to be
non-negative this implies also a bound on $|u^\eps|$. 
\epf

For the gradient estimate we aim to apply a maximum principle for
$|\nabla u^\eps|$. To do this let $f:\Omega \ra \R$ be a smooth
function and consider
$$ M=\text{graph}(f) $$
as a hypersurface in $N\times \R$, where $N\times \R$ is equipped with the
metric $\tilde{g}=\bar{g} \otimes dz^2$. Let $\nu$ be the upward pointing unit
normal of $M$ and let us take $\tau=\frac{\p}{\p z}$ as the unit
vector pointing into the upward $\R$-direction. Then define
$$ v=\big(\tilde{g}(\nu,\tau)\big)^{-1}\ .$$
As in the Euclidean case $v((p,f(p)))=\sqrt{1+|\nabla f(p)|^2}$ for 
all $p\in \Omega$. We now compute $\Delta^Mv$ at a
point $q\in M$. Let $e_1,\ldots,e_{n+1}$ be a local framing
of $TM$ around $q$ which is orthonormal at $q$. We can furthermore
assume that $\nabla^M_ve_i=0$ for all $v\in T_qM$ and $i=1,\ldots ,
n+1$. Then at $q$ we have
\beqs
\nabla^Mv = - v^2  \tilde{g}(\bar{\nabla}_{e_i}\nu, \tau) e_i = -v^2
h_{ij} \tilde{g}(e_j,\tau)e_i\ ,
\eeqs
where assume summation on $i,j$. Thus
\beq\lb{er.7}\begin{split}
\Delta^M v =&\  \text{div}^M(\nabla^M v) = \tilde{g}(\nabla^M_{e_k}\nabla^Mv,
e_k)\\
=&\ \tilde{g}\big((-2v\frac{\p v}{\p e_k}h_{ij} \tilde{g}(e_j,\tau) - 
v^2\nabla^M_{e_k}h_{ij}\tilde{g}(e_j,\tau)\\
&\ \ \ \ \  +v^2h_{ij}h_{kj}\tilde{g}(\nu,\tau))e_i,e_k\big)\\
=&\  \frac{2}{v}|\nabla^M v|^2 - v^2 \tilde{g}(\nabla^M H,\tau) - v^2
\widetilde{\text{R}}\text{ic}_{\nu k}\tilde{g}(e_k,\tau) + v|A|^2\ ,
\end{split}\eeq
where we used the Codazzi equations from the second to the third
line. To simplify further we note that
$\widetilde{\text{R}}\text{ic}(\nu,e_k) = \text{
 Ric}_N(\text{pr}_{T_qN}(\nu),\text{pr}_{T_qN}(e_k))$. We can further assume that $e_1,\ldots,e_n \perp \tau$. Then take
$$\gamma:=\frac{\text{pr}_{T_qN}(\nu)}{|\text{pr}_{T_qN}(\nu)|}\ ,$$
which is well-defined if $\nu\neq\tau$. Let us for the moment assume
that $\nu\neq\tau$. Thus
$$\text{pr}_{T_qN}(\nu) = \sqrt{1-1/{v^2}}\ \gamma\ ,\qquad
\text{pr}_{T_qN}(e_{n+1})= \pm 
\frac{1}{v}\,\gamma\ ,$$
and
\begin{equation*}
\widetilde{\text{R}}\text{ic}_{\nu k}\tilde{g}(e_k,\tau)=
\widetilde{\text{R}}\text{ic}(\nu, e_{n+1})\tilde{g}(e_{n+1},\tau)= 
-\frac{1}{v}\big(1-\frac{1}{v^2}\big)\text{Ric}_N(\gamma,\gamma).
\end{equation*}
This expression vanishes if  $\nu=\tau$, which is the right value of
the expression in (\ref{er.7}). Putting everything together we arrive at
\beq\lb{er.8}
\Delta^M v = \frac{2}{v}|\nabla^M v|^2 - v^2 \tilde{g}(\nabla^M H,\tau)
+ v |A|^2+
v\Big(1-\frac{1}{v^2}\Big)\text{Ric}_N(\gamma,\gamma)\ .
\eeq 

\bl\lb{er.l2} For any smooth solution $u^\eps$ of $(\star\star)_\eps$ the following
gradient estimate holds.
\beqs
\sup_{\Omega}|\nabla u^\eps|\leq \exp(k\,C_R\sup_\Omega
u^\eps)\cdot\sup_{\p\Omega}\big(1+ \sqrt{\eps^2+|\nabla
  u^\eps|^2}\big)\ ,
\eeqs
where $C_R:=-\inf\{\text{Ric}_N(\zeta,\zeta)\ |\ \zeta\in T_pN,\ 
  |\zeta|=1,\ p\in \Omega\}$.
\el
\bpf Examining $(\ref{er.8})$ we see that we can hope to use the
maximum principle for $v=\sqrt{1+|\nabla u^\eps|}$ if we can somehow
control the last term on the RHS. Recall that for a smooth solution 
of $(\star\star)$ the gradient bound
corresponds to a positive lower bound of the mean curvature of the
level sets $\{u=t\}= M_t$. Computing the evolution equation of
$\phi(x,t)=\exp(\eta t)H(x,t)$ for $\eta = -\min_\Omega(\text{Ric}_{N}(\nu,\nu))$
we see that 
$$ H_\text{min}(t)\geq \min(H_\text{min}(0),1)\exp(-\eta t)\ ,$$
which implies a gradient bound, given an a-priori height bound. Following
this idea we compute $\Delta^M(wv)$ where
$$w((p,z)) = \exp(-\eta z)$$
on $\Omega\times \R$ and $\eta>0$ to be chosen later. A direct computation gives
\beq\lb{er.10}
\Delta^S w = \eta^2 \Big(1-\frac{1}{v^2}\Big)w +\eta \frac{H}{v}w\
,\qquad \nabla^Sw = - \eta w \Big(\tau-\frac{1}{v}\nu\Big) \ .
\eeq
Combining this with (\ref{er.8}):
\beq\lb{er.11}\begin{split}
\Delta^M(wv) =&\ w\Delta^Mv+v\Delta^Mw +
\frac{2}{v}\tilde{g}\big(\nabla^Mv,\nabla^M(wv)\big) -
\frac{2w}{v}|\nabla^Mv|^2\\
=&\ \frac{2}{v}\tilde{g}\big(\nabla^Mv,\nabla^M(wv)\big)+ wv\Big(|A|^2
+\Big(1-\frac{1}{v^2}\Big)\text{Ric}_N(\gamma,\gamma)\\
&\  +\eta^2 \Big(1-\frac{1}{v^2}\Big)+\eta \frac{H}{v} - v
\tilde{g}\big(\nabla^MH, \tau\big)\Big).
\end{split}\eeq
Now define 
$$C_1:=\sup_{\p\Omega}\sqrt{\eps^2+|\nabla u^\eps|} $$
and assume that 
\beq\lb{er.12}
\sup_\Omega\Big(\exp\big(-\frac{\eta}{\eps} u^\eps\big)\sqrt{\eps^2 +|\nabla u^\eps|^2}\Big) >
\max\{C_1,1\} \,
\eeq
which has to be attained at an interior point. Let us take
$$ M= \text{graph}\Big(\frac{u^\eps}{\eps}\Big) \ .$$
Note that equation $(\star\star)_\eps$ implies that
\beq\lb{er.13}
H = \frac{1}{\eps^\frac{1}{k}v^\frac{1}{k}}\ ,
\eeq
where $H$ is the mean curvature of $M$. Now (\ref{er.12}) implies that
$wv$ attains an interior maximum at point $p_0$ on $S$, which is
strictly bigger than $\max\{C_1,1\}/\eps$. Furthermore by
(\ref{er.10}) and(\ref{er.13})
$$ -wv^2\tilde{g}\big(\nabla^MH,\tau\big) = \frac{1}{k}
\eps^{-\frac{1}{k}}
v^{1-\frac{1}{k}}\tilde{g}\big(\nabla^M(wv),\tau\big) +
\frac{1}{k}\eps^{-\frac{1}{k}}\eta\, w
v^{2-\frac{1}{k}}\Big(1-\frac{1}{v^2}\Big)\ .
$$
Thus at $p_0$ we have by (\ref{er.11})
\beq\lb{er.14}\begin{split}
0\geq&\  |A|^2
+\Big(1-\frac{1}{v^2}\Big)\text{Ric}_N(\gamma,\gamma)
  +\eta^2 \Big(1-\frac{1}{v^2}\Big)+
\eps^{-\frac{1}{k}}\eta v^{-1-\frac{1}{k}} \\
&\ +\frac{1}{k}\eps^{-\frac{1}{k}}\eta\,
v^{1-\frac{1}{k}}\Big(1-\frac{1}{v^2}\Big)\\
\geq&\ 
\Big(1-\frac{1}{v^2}\Big)\Big(\text{Ric}_N(\gamma,\gamma))+
\frac{\eta}{k\eps}(\eps v)^{1-\frac{1}{k}}\Big) \ .
\end{split}
\eeq
If we now choose $ \frac{\eta}{\eps}= k \cdot C_R$
we arrive at a contradiction since at $p_0$ we have $ \eps v >  w^{-1}\max\{C_1,1\}\geq 1$.
\epf

\bl\lb{er.l3} Let $\min_{\p\Omega} H_{\p\Omega} := \delta_0>0$ and
$0<\delta_1\leq \delta_0/(2C_R)$ be such that $d(p):=\text{dist}(p,\Omega)$ is smooth on 
$\Omega_{\delta_1}=\{p\in\Omega\ | \ d(p)<\delta_1\}$. 
Let $0<\eps<\eps_0$ where $\eps_0:=\min\{C_2,1\}$, $C_2:= \sup\{2^{\frac{3k-1}{2}}\delta_0^{-k},
\delta_1^{-1}C_1\}$ and $C_1$ is the a-priori bound on $\sup_\Omega
|u^\eps|$ from Lemma \ref{er.l1}. Then any smooth solution $u^\eps$
of $(\star\star)_\eps$ satisfies the estimate
\beq\lb{er.15}
\sup_{\p\Omega}|\nabla u^\eps| \leq C_2\ .
\eeq
\el 
\bpf We construct a barrier at the boundary, and since $u^\eps \geq 0$ we
only need a barrier from above. To construct a suitable supersolution $\Phi$,
we try the ansatz $\Phi(p) = \beta\cdot d(p)$ for a constant
$\beta>0$. Observe that on $\Omega_{\delta_1}$ we have $\Delta d =
-H_{S_r}$, where $H_{S_r}$ is the mean curvature of the hypersurfaces
$S_r:=\{d=r\}$. Computing as in the proof Lemma \ref{er.l1} we find that
a sufficient condition for $\Phi$ to be a supersolution
$\Omega_\delta$ is that
\beq\lb{er.16}
H_{S_r} \geq \frac{1}{\beta}(\eps^2+\beta^2)^\frac{k-1}{2k} \ ,
\eeq
for all $0<r<\delta$ and a suitable $0<\delta<\delta_1 $. The evolution
equation of $H_{S_r}$ along a geodesic is given by
$$\frac{\p}{\p r}H_{S_r}= |A_{S_r}|^2 + \text{Ric}(\nu,\nu) \geq -C_R\ ,$$
which implies
$$H_{S_r}\geq H_{\p\Omega}-C_Rr \geq \delta_0-C_Rr\geq
\frac{\delta_0}{2}$$
for $0\leq r \leq \delta_1$. Assuming that $\eps\leq \beta$, a
sufficient condition to fulfill (\ref{er.16}) is that $\beta \geq
2^\frac{3k-1}{2}\delta_0^{-k}$. If we further assume that $\beta \geq
C_1/\delta_1$ we can ensure that
$$ \Phi \geq u^\eps$$
on $S_{\delta_1}$. Thus by the maximum principle $\Phi \geq u^\eps$ on
$\Omega_{\delta_1}$, which gives the desired gradient estimate.
\epf

To show the existence of solutions to $(\star\star)_\eps$ we study
solutions to the following family of equations
\begin{equation*}
 (\star\star)_{\eps,\kappa} \begin{cases}
 \ \ \text{div}_{N}\Big(\frac{\nabla u^{\eps,\kappa}}
{\sqrt{\eps^2+|\nabla u^{\eps,\kappa}|^2}}
\Big ) =-\kappa\big(\eps^2+|\nabla u^{\eps,\kappa}|^2\big)^{-\frac{1}{2k}}\ \ 
&\text{in} \ \Omega\\[+2.0ex]
\ \ u^{\eps,\kappa} = 0 &\text{on}\ \p\Omega
\end{cases}\qquad\qquad
\end{equation*}
for $0\leq \kappa \leq 1$ and $0<\eps<\eps_0$. In the following we
show that for any fixed $0<\eps<\eps_0$ we have uniform a-priori
sup and gradient estimates in $\kappa$. Since $\kappa\leq 1$ it is
easy to check that $(\ref{er.2})$ and $(\ref{er.15})$ also hold for 
smooth solutions of $(\star\star)_{\eps,\kappa}$: 
\beq\lb{er.16a}
\sup_{\Omega}|u^{\eps,\kappa}| \leq C(n,k,\text{diam}(\Omega))\ ,
\qquad
\sup_{\p\Omega}|\nabla u^{\eps,\kappa}|\leq C_2 \ ,
\eeq
for all $0\leq \kappa\leq 1$. Here $C_2$ is the constant from Lemma
\ref{er.l3}, where we assume the same conditions on $\p\Omega$. 
For the interior gradient estimate, fix an $\eps \in
(0,\eps_0)$. Let us work on the hypersurface
$$ M=\text{graph}\Big(\frac{u^{\eps,\kappa}}{\eps}\Big) \ .$$
Equation $(\star\star)_{\eps,\kappa}$ then implies that the mean
curvature $H$ of $S$ is given by
$$ H = \frac{\kappa}{\eps^\frac{1}{k}v^\frac{1}{k}}\ .$$
As in the proof of Lemma \ref{er.l1} we obtain that at an interior
maximum of $v$, we have the inequality (compare $(\ref{er.14})$)
\beqs
\begin{split}
0\geq&\  |A|^2
+\Big(1-\frac{1}{v^2}\Big)\text{Ric}_N(\gamma,\gamma)
  +\eta^2 \Big(1-\frac{1}{v^2}\Big)+
\kappa\,\eps^{-\frac{1}{k}}\eta v^{-1-\frac{1}{k}} \\
&\ +\kappa\, \frac{1}{k}\eps^{-\frac{1}{k}}\eta\,
v^{1-\frac{1}{k}}\Big(1-\frac{1}{v^2}\Big)\\
\geq&\ 
\Big(1-\frac{1}{v^2}\Big)\big(\text{Ric}_N(\gamma,\gamma))+ \eta^2\big) \ ,
\end{split}
\eeqs
which gives a contradiction if $\eta > \sqrt{C_R}$ and $v>1$. This
yields the interior estimate 
\beq\lb{er.18}
\sup_\Omega |\nabla u^{\eps,\kappa}| \leq
\exp\big(\sqrt{C_R}\,\eps^{-1}\sup_\Omega u^{\eps,\tau}\big)\cdot
  \sup_{\p\Omega}\big(\sqrt{\eps^2+|\nabla u^{\eps,\tau}|^2}\big) \ ,
\eeq
for all $0\leq\kappa\leq 1$.

\bl\lb{er.l4}Under the assumptions of Lemma \ref{er.l3} a smooth
solution to $(\star\star)_\eps$ exists. 
\el

\bpf We aim to apply the method of continuity to
$(\star\star)_{\eps,\kappa},\ 0\leq \kappa\leq 1$. Fix an $\eps \in
(0,\eps_0)$ and write $(\star\star)_{\eps,\kappa}$ as
$$ F^\kappa(w):=\text{div}_N\Big(\frac{\nabla
  w}{\sqrt{\eps^2+|\nabla w|^2}}\Big) +\kappa \big(\eps^2 +|\nabla
w|^2 \big)^{-\frac{1}{2k}}=0\ , $$
with $w=0$ on $\p \Omega$. The map 
$$ F\, :\, C_0^{2,\alpha}(\bar{\Omega})\times \R \ra
C^\alpha(\bar{\Omega})\ ,$$
defined by $F(w,\kappa):= F^\kappa(w)$ is $C^1$ and possesses the
solution $F(0,0)=0$. Let 
$$ I:=\{\kappa \in [0,1]\, |\, (\star\star)_{\eps,\kappa}\ \text{has a
solution}\ u^{\eps,\kappa}\in C^{2,\alpha}_0(\bar{\Omega})\} \ .$$
Clearly $0\in I$. We want to show that $I$ is relatively open and
closed.\\[+1ex]
To see that $I$ is closed we note that the a-priori estimates
$(\ref{er.16a})$ and $(\ref{er.18})$ imply uniform
$C^1(\bar{\Omega})$-bounds 
and thus $(\star\star)_{\eps,\kappa}$ is uniformly elliptic. The
Nash-Moser-DeGiorgi estimates then yield uniform bounds in
$C^{1,\alpha}(\bar{\Omega})$. Applying Schauder estimates we obtain
uniform bounds in $C^{k,\alpha}(\bar{\Omega})$ for any $k\geq 2$. Thus
by the Arzela-Ascoli theorem $I$ is closed.\\[+1ex]
To prove that $I$ is also open we linearize the map $F^\kappa$ at a
solution $u$. This linearization is given by 
$$DF^\kappa|_u\, : \, C^{2,\alpha}_0(\bar{\Omega})\ra
C^{\alpha}(\bar{\Omega})\ . $$
Now $F^\kappa(w)$ has the form
$$F^{\kappa}(w)= \nabla_iA^i(\nabla w) + B(\nabla w)\ ,$$
which is independent of $w$. So by the maximum principle the
linearization
$$DF^\kappa|_u(v)= \nabla_i(A^i_{p_j}(\nabla u)\nabla_jv) +
B_{p_j}(\nabla u) \nabla_jv $$
has only the zero solution. Using linear existence theory and Schauder
estimates up to the boundary we see that $DF^\kappa|_u$ is an
isomorphism.
So by the implicit function theorem, the set of $\kappa$ such that
$F(u,\kappa)=0$ has a solution (namely $I$) is open. Therefore $1\in
I$ which proves the existence of $u^\eps$ in $C^{2,\alpha}$. Smoothness
follows again by Schauder estimates.
\epf


\section{Weak $H^k$-flow and Uniqueness}\lb{uniqueness}

\noindent Given a connected, open and bounded set $\Omega\subset N$ with smooth
boundary, s.t. $H|_{\p\Omega}>0 $, the results of the last section
ensure the existence of smooth solutions $u^\eps$ to $(\star\star)_\eps$ for
sufficiently small $\eps>0$. The a-priori estimates guarantee uniform
bounds in $C^1(\bar{\Omega})$, independent of $\eps$. Thus, given any
sequence $\eps_i \ra 0$, we can extract a subsequence (again denoted by
$(\eps_i)$), such that
$$u^{\eps_i}\ra u $$
in $C^0(\bar{\Omega})$ to a function $u\in
C^{0,1}(\bar{\Omega})$. This suggests the following definition:

\bd Let $\eps_i \ra 0$ and corresponding solutions $u^{\eps_i}$ to
$(\star\star)_{\eps_i}$ be given. Assume that $u^{\eps_i}\ra u$
uniformly on $\bar{\Omega}$, where $u, u^{\eps_i}$ are uniformly
bounded in $C^{0,1}(\bar{\Omega})$. We
then call $u$ a weak $H^k$-flow with initial condition $\Omega$.
\ed

By the reasoning above we have the existence of a weak $H^k$-flow. We now
want to show that such a weak solution actually is unique, i.e. the
limiting function $u$ is independent of the approximating sequence $u^{\eps_i}$. To do this
we first want to show that any weak $H^k$-flow coincides with the
smooth flow, as long as the latter exists. 

Let $F(\cdot,t):\p\Omega\times [0,T)\ra N$ be the unique
solution to $(\star)$ with initial condition
$F(\cdot,0)=\text{Id}_{\p\Omega\ra \p\Omega}$. We may further assume
that $T>0$ is maximal. Comparing with shrinking distance-spheres as in
the proof of Lemma \ref{er.l1} we see that $T$ is finite. 
Let us write $M_t=F(\p\Omega,t)$. As
in the beginning of the proof of Lemma \ref{er.l2} we furthermore
obtain that the mean
curvature of the hypersurfaces remains strictly positive. Thus we 
have that $M_{t_1}\cap M_{t_2} = \emptyset$ if $t_1 \neq
t_2$. For $0<\tau\leq T$ define
$$\Omega_\tau=\bigcup_{0<t<\tau}\!\!\!M_t \subset \Omega\ $$
and $u^*:\Omega_\tau \ra \R^+$ by $u^*(p)=t$, if $x \in M_t$.

\bl\lb{eu.l1} Let $u$ be a weak solution with initial condition
$\Omega$. Then
$$u=u^*$$
in $\Omega_T$.
\el
\bpf Choose any $0<\tau<T$. Then $u^*$ is a smooth solution of
$(\star\star)$ on $\overline{\Omega}_\tau$. Our aim is to twist $u^*$
somewhat to construct upper and lower barriers for
$(\star\star)_\eps$ for $\eps$ small enough. To do this let $\psi:[0,\tau) \ra \R^+$ be a smooth,
increasing function with $\psi(0)=0$ and define
$$ v(x) = \psi(u^*(x)) $$
for $x \in \Omega_\tau$. Since $u^*$ is a solution of $(\star\star)$ we
see by a direct calculation that in $\Omega_\tau$
\beq\lb{eu.1}
\text{div}_N\bigg(\frac{\nabla v}{|\nabla v|}\bigg) =
\frac{1}{|\nabla v|}\bigg(\delta^{ij}-\frac{\nabla^iv \nabla^j
  v}{|\nabla v|^2}\bigg) \nabla_i\nabla_j v = -
\frac{(\psi^\prime)^\frac{1}{k}}{|\nabla v|^\frac{1}{k}}\ ,
\eeq
where in the second equality we use normal coordinates at $p\in
\Omega_\tau$. This implies that $v$ is a subsolution to $(\star\star)$
if $\psi^\prime \leq 1$ and a supersolution if $\psi^\prime \geq
1$. Using (\ref{eu.1}) we furthermore have in normal coordinates at a
point $p\in \Omega_\tau$ that
\beq\lb{eu.2}\begin{split}
\bigg(\delta^{ij}&\ -\ \frac{\nabla^iv\nabla^jv}{\eps^2+|\nabla
    v|^2}\bigg)\nabla_i\nabla_j v + \Big(\eps^2 +|\nabla
v|^2\Big)^\frac{k-1}{2k} =\\
&\ \Big(\big(\eps^2+|\nabla v|^2\big)^\frac{k-1}{2k}-|\nabla
v|^\frac{k-1}{k}\Big)
 + \bigg(\frac{\nabla^iv\nabla^j v}{|\nabla v|^2} -
\frac{\nabla^iv\nabla^jv}{\eps^2 +|\nabla v|^2}\bigg)\nabla_i\nabla_j
v \\
&\ + \big(1-(\psi^\prime)^\frac{1}{k}\big) |\nabla v|^\frac{k-1}{k}\ .
\end{split}
\eeq
Now on $\Omega_\tau$ there exists a positive constant $C_1<\infty$ such
that 
\beqs
\frac{1}{C_1} \leq |\nabla u^*| \leq C_1\ .
\eeqs
The positive lower bound on the gradient is due to the fact that the
flow $F(\cdot, t)$ is smooth up to $\tau$. The upper bound comes from
the uniform positive lower bound on the mean curvature of the surfaces
$M_t$.\\[0.5ex]
We first want to construct an upper barrier for
$(\star\star)_\eps$. Pick any sufficiently small $\delta>0$ and take
$\psi(r)=(1+\delta)r$ for $r\in [0,\tau-\delta]$ and continue $\psi$
smoothly in a convex way on $(\tau-\delta,\tau]$, such that 
\beq\lb{eu.3}\frac{\psi^\prime(\tau)}{C_1}\geq \sup_{0\leq \eps\leq
  \eps_0}\sup_\Omega |\nabla u^\eps| + 1 
\eeq
where $\eps_0$ is the constant from Lemma \ref{er.l3}. Observe that 
$|\nabla v|= \psi^\prime |\nabla u^*|$. Thus (\ref{eu.3}) implies that
the solutions $u^\eps$ cannot touch $v$ in a neighborhood of the inner
boundary of $\Omega_\tau$. Having fixed $\psi$, we have the bounds
$$\frac{1}{C_2}\leq |\nabla v| \leq C_2\ , \quad |\nabla^2v|\leq C_3\
\text{on}\ \Omega_\tau $$
for some positive constants $C_2,C_3$. Since $\psi^\prime\geq
1+\delta$, we see from equation (\ref{eu.2}) that $v$ is a
supersolution on $\Omega_\tau$ of $(\star\star)_\eps$ for
sufficiently small $\eps$. As explained above, the solutions
$u^{\eps_i}$ which converge to $u$ cannot touch $v$ on the inner
boundary of $\Omega_\tau$, so $v$ acts as an upper barrier for
sufficiently small $\eps_i$. Taking the limit $u^{\eps_i}\ra u$ we obtain that
$u\leq v$ on $\Omega_\tau$. Letting $\delta \searrow 0$ and $\tau \nearrow
T$ we arrive at $u\leq u^*$ on $\Omega_T$. \\[1.5ex]
For the lower barrier take again $\tau, \delta$ as above and let 
$\psi(r)=(1-\delta)r$ for $r\in [0, \tau
-2\delta]$. The aim is then to continue $\psi$ on $(\tau
-2\delta,\tau]$ in a concave way, such that we can extend $v$ by a
constant on $\Omega\setminus \Omega_\tau$ to obtain a
$C^2$-subsolution to $(\star\star)_\eps$ on the whole of
$\Omega$. To estimate the RHS of (\ref{eu.2}) from below we drop the
first term and use that $\nabla v = \psi^\prime \nabla u^*$,\ 
$\nabla_i\nabla_j v = \psi^{\prime\prime}\nabla_iu^*\nabla_ju^*+
\psi^\prime \nabla_i\nabla_j u^*$. Thus we get a lower estimate of
the RHS by
\beq\lb{eu.4}\begin{split}
\bigg(1&-\frac{(\psi^\prime)^2|\nabla
  u^*|^2}{\eps^2+(\psi^\prime)^2|\nabla
  u^*|^2}\bigg)\Big(\psi^{\prime\prime}|\nabla u^*|^2 - \psi^\prime
|\nabla^2 u^*|\Big)\\ &\ +\ 
(\psi^\prime)^\frac{k-1}{k}\big((1-(\psi^\prime)^\frac{1}{k})|\nabla
u^*|^\frac{k-1}{k}\big)
\geq C_1^2\psi^{\prime\prime}-C_3\psi^\prime + \gamma(\psi^\prime)^{\frac{k-1}{k}}\ ,
\end{split}\eeq
where $\gamma >0$ is a constant depending only on $\delta$ and
$C_1$. Thus the RHS is nonnegative on $[\tau-\delta,\tau]$, if 
\beqs
\psi^\prime\leq \Big(\frac{\gamma}{2C_3}\Big)^{k}\ \
\text{and}\ \ 
-\psi^{\prime\prime} \leq \frac{1}{2C_1^2}\gamma 
(\psi^\prime)^\frac{k-1}{k}
\eeqs
It is then a direct calculation that the choice $\psi(r)= -\alpha
(\tau-r)^{k+1}+b$ satisfies the above constraints for a constant
$\alpha>0$, depending on $k,\delta,\gamma, C_1,C_3$, and $b$ still free
to choose. We then adjust $b$ such that we can continue $\psi$ on
$[\tau-2\delta,\tau-\delta]$ smoothly in a concave way. On
$\Omega_{\tau-\delta}$ we can again use (\ref{eu.2}) to see that for
small enough $\eps$ the function $v$ is a subsolution. Equation
(\ref{eu.4}) then guarantees that this also works on
$\Omega_\tau\setminus\Omega_{\tau-\delta}$. We
extend $v$ to the whole of $\Omega$ by setting $v(p)=b$ for
$p\in\Omega\setminus \Omega_\tau$. Thus for $k>1,\ v$ is a $C^2$
subsolution of $(\star\star)_\eps$ for sufficiently small $\eps$ and
thus a lower barrier for $u^{\eps_i}$. Using the a-priori estimates we
can take the limit $k\searrow 1$ to see that $v$ also acts as a lower
barrier in this case. Then arguing as in the case of the upper barrier
we obtain finally that $ u\geq u^*$ on $\Omega_T$.
\epf

\bc[Avoidance of smooth flows]\lb{eu.c1} Let $u$ be a weak $H^k$-flow with initial
condition $\Omega$ and $(M_t)_{t_0\leq t \leq t_1},\ t_0\geq 0,$ be a smooth, compact
$H^k$-flow with positive mean curvature. Assume that $M_{t_0}$ and
$u$ are disjoint at $t_0$, i.e. $M_{t_0}\cap \{u=t_0\} = \emptyset$,
then they remain so for all future times, i.e. $M_{t}\cap \{u=t\} =
\emptyset,\ \forall\ t_0\leq t\leq t_1 $.
\ec
\bpf Let $\Omega^\prime$ be the bounded and open set in $\rne$ such
that $M_{t_0}=\p\Omega^\prime$. We can assume that $\Omega^\prime
\subset \Omega$ or $\Omega\subset \Omega^\prime$, otherwise there 
is nothing to prove. We treat the case that $\Omega^\prime \subset
\Omega$, the other case follows similarly. Let $u^{\eps_i}$
be the sequence of solutions to $(\star\star)_{\eps_i}$ converging to
$u$. Then take $\tilde{u}^{\eps_i}$ to be the solutions to
$(\star\star)_{\eps_i}$ on $\Omega^\prime$, where we take the same
sequence $\{\eps_i\}$. We can furthermore assume that
$\tilde{u}^{\eps_i}\ra \tilde{u}$ uniformly, s.t. $\tilde{u}$ is a
weak solution with initial condition $\Omega^\prime$. By Lemma
\ref{eu.l1} we have that 
\beq\lb{eu.6}
M_t = \{\tilde{u}=(t-t_0)\}
\eeq
for $t_0\leq t \leq t_1$. Since $u > t_0$ on $\Omega^\prime$ 
also $u^{\eps_i}> t_0$
for sufficiently large $i$, and thus by the maximum principle
$u^{\eps_i}\geq \tilde{u}^{\eps_i}+t_0$ which gives
$$ u\geq \tilde{u}+t_0$$
on $\Omega^\prime$. We can now shift $M_{t_0}$ a little bit forward in
time such that $M_{t_0}$ remains disjoint from $u$
and repeat the above argument. This yields
$$ u > \tilde{u}+t_0$$
on $\Omega^\prime$, which proves the claim, using (\ref{eu.6}).
\epf

By interposing a $C^{1,1}$-hypersurface between two weak $H^k$-flows, 
we can argue as it is done in \cite{Ilmanen93} for set-theoretic
subsolutions to mean curvature
flow, that also two weak $H^k$-flows satisfy the avoidance
principle. Since the proof depends on the translational invariance of
the flow it works only if the surrounding space is Euclidean. Furthermore
we can start a smooth $H^k$-flow only from a
$C^{1,1}$-hypersurface if it has nonnegative mean curvature, so we have
to restrict ourselves to low dimensions.

\bt\lb{eu.t1} Let $N=\rne,\ n\leq 6,$ and $u,\tilde{u}$ be two weak $H^k$-flows generated by
two open sets $\Omega,\tilde{\Omega}\subset \rne$ where at least one
of them is bounded. Assume that for
$t_1,t_2\geq 0$ we have $\{u=t_1\}\cap \{\tilde{u}=t_2\}
=\emptyset$. Then $\{u=(t_1+\tau)\}\cap \{\tilde{u}=(t_2+\tau)\}
=\emptyset$ for all $\tau\geq 0$. 
\et

\bpf We can assume that $\{\tilde{u}\geq t_2\} \Subset \{u>t_1\}$. We
want to show that
\mbox{$\rm{dist}\big(\{u=(t_1+\tau)\},\{\tilde{u}=(t_2+\tau)\}\big)$} is
non-decreasing in $\tau$. Observe that by translational invariance of
the $H^k$-flow Corollary \ref{eu.c1} proves this if one of the two
flows is actually smooth. By \cite{Ilmanen93} we know that there exists
a closed $C^{1,1}$-hypersurface $S\subset \{u>t_1\}\setminus \{\tilde{u}\geq t_2\}$
which separates $\{u=t_1\}$ and $\{\tilde{u}=t_2\}$ such that
\beq\lb{eu.7}\text{dist}\big(\{u=t_1\},\{\tilde{u}=t_2\}\big)=
\text{dist}\big(\{u=t_1\},
S\big)+\text{dist}\big(\{\tilde{u}=t_2\}, S\big)\ .
\eeq
Let $E$ be the open set, bounded by $S$, such that
$\{\tilde{u}\geq t_2\}\subset E$. Now take $E^\prime \subset
\{u>t_1\}$ to
be the outer minimizing hull of $E$ in $ \{u>t_1\}$, i.e. $E^\prime$ is the
intersection of all minimizing hulls in $ \{u>t_1\}$, which contain
$E$. Here we call a set $F\subset G$, where $G\subset \rne$ is
open, a minimizing hull in $G$, if it minimizes area from the outside in
$G$, that is, if
\beq\lb{eu.8}|\p^*F\cap K| \leq |\p^*H\cap K|
\eeq
for any $H$ containing $F$ such that $H\setminus F \Subset G$, and any
compact set $K$ containing $H \setminus F$. For details on minimizing
hulls, see \cite{HuiskenIlmanen}, chapter 1. Since $u$ is weak $H^k$-flow,
all the sets $\{u > t\}$ are minimizing hulls in $\Omega$, see
Corollary \ref{pr.c2}, which is proved
independently of this uniqueness result. By Corollary \ref{eu.c1} the
sets $\{u=t\}$ cannot develop an interiour, thus there is a $\tau>t_1$
such that $E \Subset\{u>\tau\}$. This implies that $E^\prime$ cannot touch
$\p\{u>t_1\}$. Since $\p E = S$ is $C^{1,1}$ and $n\leq 6$, a result of Sternberg,
Williams and Ziemer \cite{SternbergWilliamsZiemer91} implies that
$S^\prime:=\p E^\prime$ is also a $C^{1,1}$-hypersurface. Since
$E^\prime$ locally minimizes area from the outside, $S^\prime$ carries
a nonnegative weak mean curvature which is in $L^\infty$. By an
argument of Huisken and Ilmanen, see \cite{HuiskenIlmanen}, Lemma 2.5, $S^\prime$
can be approximated by a sequence $S^\prime_i$ of smooth
hypersurfaces from the inside, which are uniformly controlled in $C^2$ and with strictly positive
mean curvature. Let $M^i_t$ be the smooth evolution along the
$H^k$-flow of the $S_i^\prime$'s. By Proposition 3.9 in \cite{convex}
these flows exist on a uniform time interval $[0,\eps)$, for some
$\eps>0$. We first want to show that
\beqs
\text{dist}\big(\{u=t_1\},S\}\big) = \text{dist}\big(\{u=t_1\},S^\prime\}\big)\ .
\eeqs 
Assume to the contrary that
$$ \text{dist}\big(\{u=t_1\},S\}\big) >
\text{dist}\big(\{u=t_1\},S^\prime\}\big)\ ,$$
then the shortest distance has to be attained at a point $p\in
S^\prime$ such that $S^\prime$ is a smooth minimal surface around $p$. Thus the $M^i_t$
move initially near $p$ as slowly as we wish, and note that $u$ has to
avoid the flows $M^i_t$. Even more the shortest distance from
$\{u=t_1\}$ to $M^i_0$ has to be attained at a sequence of points $p^i\in
M^i_0$, which we can assume to converge to $p$. Translating the $M^i_0$'s such that they touch
$\{u=t_1\}$ in $p^i$, we get by the avoidance principle w.r.t.\!\! smooth
flows a contradiction to the $C^{0,1}$-bound of
$u$. Thus we can replace $S$ by $S^\prime$ in (\ref{eu.7}).
Now $u$ and $\tilde{u}$ avoid the evolution of all the $M^i_t$ which
proves the statement of the theorem in the limit $i\ra \infty$.
\epf

By shifting the initial condition a little bit in time this implies
uniqueness.

\bc\lb{eu.c2} Let $N=\rne,\ n\leq 6$ and $\Omega \subset
\rne$ be a bounded, open set with smooth boundary such that
$H_{\p\Omega}>0$. Then the weak $H^k$-flow generated by $\Omega$ is unique.
\ec

Note that in the case $k=1$, i.e.\! mean curvature flow, any
weak flow as above, without any restriction on the dimension and
on the ambient space, coincides with the level-set flow of $\p\Omega$
and is thus unique, see \cite{EvansSpruckI}. 

\section{Further Properties and Regularity}\lb{regularity}

\noindent In the following section let $\Omega \subset N$ be a fixed open
and bounded set with smooth boundary such that $H_{\p\Omega}>0$. Let
$u\in C^{0,1}(\overline{\Omega};\R^+)$ be a weak $H^k$-flow generated by
$\Omega$, i.e. there exists a sequence $\eps_i \searrow 0$ and
solutions $u^{\eps_i}$ to $(\star\star)_{\eps_i}$ which are uniformly
bounded in $C^{0,1}(\overline{\Omega};\R^+)$ converging to $u$ in
$C^0(\overline{\Omega})$. Then the hypersurfaces \mbox{$N^i_t\subset
N\times \R$}, defined by 
$$
N^i_t:=N^{\eps_i}_t=\text{graph}\bigg(\frac{u^{\eps_i}}{\eps_i}-\frac{t}{\eps_i}\bigg)\
,$$
which are level sets $\{U^{\eps_i}=t\}$ of the function
\mbox{$U^{\eps_i}\big((x,z)\big)=u^{\eps_i}(x)-\eps_i z$} on $\Omega\times
\R$, are smooth translating solutions of the  $H^k$-flow $(\star)$, see
(\ref{er.1}). Equation $(\star\star)_{\eps_i}$ implies
that the mean curvature $H^i_t$ of $N^i_t$ is given by
\beq\lb{pr.1}
H^i_t =  \frac{1}{\big( \eps_i^2+|\nabla
  u^{\eps_i}|^2\big)^\frac{1}{2k}}\ .
\eeq
To fix some further notation define the following subsets of
$\Omega\times \R$:
\beqs
E^i_t:=\{U^{\eps_i}>t\},\qquad E^\prime_t:= \{U>t\}\ ,
\eeqs
where $U\big((x,z)\big)=u(x)$ on $\Omega\times \R$. The sets
$E^\prime_t$ can be written as $E^\prime_t=E_t\times \R$, where
$E_t:=\{u>t\} \subset \Omega$. A first observation is that the sets
$E^i_t$ are minimizing hulls in $\Omega\times \R$, see (\ref{eu.8}).

\bl\lb{pr.l1} The sets $E^i_t$ are minimizing area from outside in
$\Omega\times \R$, that is 
$$|\p^*E^i_t\cap K|\leq |\p^*F\cap K|$$
for $F$ with $E^i_t\subset F,\ F\setminus E^i_t \subset K \subset
\Omega\times \R$, where $K$ is compact. Here we take the right hand
side to be $+\infty$ if $F$ is not a Caccioppoli-set. 
\el

\bpf The outward unit normal to the surfaces $N^i_t$, which is given
by $\nu= -\nabla U^{\eps_i}/|\nabla U^{\eps_i}|$ is a smooth
vectorfield on $\Omega \times \R$ with
div$(\nu)=|DU^{\eps_i}|^{-1/k}>0$. Thus we get by the divergence
theorem for Caccioppoli-sets, using $\nu$ as a calibration:
\begin{equation}\lb{pr.3}\begin{split} 0 \leq \!\! \inl_{F\setminus E^i_t}\!\!\!\! \text{div}(\nu)\,
dx =&\  - \!\!\!\!\! \inl_{
\p^*E^i_t\cap K}\!\!\!\!\!\nu \cdot \nu_{\p^*E^i_t}\, d\hn + \!\!\! \inl_{\p^*F \cap K}
\!\!\!\!\!\nu \cdot \nu_{\p^*F}\, d\hn \\
\leq&\  -|\p^*E^i_t\cap K|\  +\  |\p^*F\cap K| ,
\end{split}
\end{equation}
where we take $\nu_{\p^*E^i_t},\ \nu_{\p^*F}$ to be the outward unit
normals to $\p^*E^i_t,\ \p^*F$.
\epf

\bc{(Mass Bound).}\lb{pr.c1} Let $\Omega^\prime:=\Omega \times
[a,b]$. Then 
\beq\lb{pr.4}|N^i_t\cap \Omega^\prime| \leq (b-a)\ \hn(\p\Omega) + 2 \hne(\Omega)\eeq
for all $-\infty < t< \infty$. 
\ec
\bpf Let $\Omega_j\Subset \Omega$ be a sequence of open sets, such
that $\p\Omega_j \ra \p\Omega$ in $C^1$. Then $F_j:= \big(\Omega_j\times
(b-a)\big)\cup E_t^i$ are valid comparison sets, and the above lemma
gives the estimate for $j \ra \infty$.
\epf 

We can use this a-priori mass bound and the lower bound on the mean
curvature together with evolution equations to deduce space-time
bounds, independent of $\eps_i$.

\bl\lb{pr.l2} Let $I = [a,b]\subset \R $ be a bounded interval. Then
\beq\lb{pr.5}
 \inl_{-\infty}^{+\infty}\!\!\int_{N^i_t \cap (\Omega \times I)}\!\!
H_i^{k+1}\, d\hne \leq (b-a)\,\hn(\p \Omega) + 2 \hne(\Omega)\ .
\eeq
\el
\bpf Observe that by the Coarea formula
$$ \inl_{\Omega\times I} |DU^{\eps_i}|^{-\frac{1}{k}} \, dx =
\inl_{-\infty}^{+\infty}\!\!\int_{N^i_t \cap (\Omega \times I)}
H_i^{k+1}\, d\hne \ .$$
We can then use (\ref{pr.3}) and argue as in in the proof of Corollary
\ref{pr.c1}.
\epf

\bl\lb{pr.l3} Let $k>3$, $I = [a,b]\subset \R $ be a bounded interval
and $\Omega^\prime \Subset \Omega$. Then 
\beqs 
\inl_{-\infty}^{+\infty}\!\!\int_{N^i_t \cap (\Omega^\prime \times I)}\!\!
|\nabla H_i|^2\, d\hne \leq  C(\Omega, \Omega^\prime, I, k)\ .
\eeqs
\el
\bpf Choose $\phi \in C^2_c(\Omega\times \R),  0\leq \phi \leq 1,$ such that $\phi = 1$ on
$\Omega^\prime\times I$ and let $\alpha >0$. By the evolution equation for
the mean curvature 
and integration by parts we can compute 
\begin{equation*}\begin{split}
\dt \inl_{N^i_t} \phi H_i^{-\alpha}\, d&\hne =\
\inl_{N^i_t}(-\alpha)(1+\alpha)k\phi H_i^{k-\alpha -3} |\nabla
H_i|^2\\
&\ + \alpha k H_i^{k-\alpha -2} \langle \nabla H_i,\nabla \phi\rangle 
 -\alpha \phi
 H_i^{k-\alpha-1}(|A|^2+\overline{\text{R}}\text{ic}(\nu,\nu))\\[1ex]
 &\ - \langle\nabla \phi, \nu\rangle
H_i^{k-\alpha} - \phi H_i^{k+1-\alpha}\, d\hne\ .
\end{split}\end{equation*}
We can estimate the first term in the second line as follows:
$$\langle \nabla H_i,\nabla \phi \rangle H_i^{k-\alpha-2} \leq
\frac{1}{2}\frac{|\nabla \phi|^2}{\phi} H_i^{k-\alpha -1} +
\frac{1}{2} \phi H_i^{k-\alpha -3}|\nabla H_i|^2\ . $$
Since $|\nabla \phi|^2/\phi \leq C(\phi)$ we can apply this estimate
and integrate from $t_1$ to $t_2$ to  arrive
at
\begin{equation*}\begin{split}
\inl_{t_1}^{t_2} k \phi H_i^{k-\alpha -3}|\nabla H_i|^2\, d\hne\, dt
\leq&\  \frac{2}{\alpha}\inl_{N^i_{t_1}}\phi H_i^{-\alpha}\, d\hne 
-  \frac{2}{\alpha}\inl_{N_{t_2}^i}\phi H_i^{-\alpha}\, d\hne \\
+&\ C\inl_{t_1}^{t_2}\!\!\int_{N^i_t} H_i^{k-\alpha-1}
+ H_i^{k - \alpha}\, d\hne\, dt\ ,
\end{split}\end{equation*}
where $C=C(\phi, \alpha^{-1}, k, \sup_\Omega|\overline{\text{R}}\text{ic}| )$.
Now for $t_1\ll -1,\ t_2 \gg 1$ we have $N^i_{t_1}\cap \text{supp}(\phi)=
N^i_{t_2}
\cap \text{supp}(\phi)= \emptyset$ such that the first two terms on
the RHS drop out. Furthermore the a-priori gradient bound from Lemmata
\ref{er.l1}, \ref{er.l2} and \ref{er.l3} together with (\ref{pr.1})
give a uniform positive lower bound on $H_i$. Now choose $\alpha = (k-3)/2$ and use Lemma
\ref{pr.l2} to prove the claimed estimate.
\epf

The space-time estimate (\ref{pr.5}) implies that the measure of the
sets $E_t$ is H\"older-continuous in time, which also excludes that
the level sets of $u$ can ``fatten up'': 

\bl\lb{pr.l4} The weak $H^k$-flow $u$ is non-fattening, i.e. $\hne
\big(\{u=t\}\big)=0$ for all $t\in [0,T]$, where $T=\sup_\Omega u$. 
\el
\bpf Let $\Omega^\prime := \Omega \times (0,1)$. Then for any $t_1,t_2
\in \R,\ t_1 <t_2$ we have by the Coarea formula together with
(\ref{pr.4}), (\ref{pr.5})  and  H\"older's inequality
\begin{align*}|\hnz(E^i_{t_1}\cap \Omega^\prime) -&\ \hnz(E^i_{t_2}\cap
\Omega^\prime)| = \inl_{t_1}^{t_2} \!\!\int_{N_t^i\cap\Omega^\prime} H_i^k\,
d\hne\, dt \\ 
\leq&\ |t_2-t_1|^{\frac{1}{k+1}}
\Bigg(\inl_{t_1}^{t_2}\bigg(\int_{N_t^i\cap \Omega^\prime} H_i^{k}\,
d\hne\bigg)^\frac{k+1}{k} dt\Bigg)^\frac{k}{k+1}\displaybreak[0]\\
\leq&\ C\, |t_2-t_1|^\frac{1}{k+1} \Bigg(\inl_{t_1}^{t_2}\int_{N_t^i\cap \Omega^\prime} H_i^{k+1}\,
d\hne\, dt\Bigg)^\frac{k}{k+1}\displaybreak[0]\\
\leq&\  C |t_2-t_1|^\frac{1}{k+1}\ ,
\end{align*}
where the constant $C$ does not depend on $i$. Observe that since
$U^{\eps_i}\ra U$ locally uniformly on $\Omega\times \R$ we have that
\beq\lb{pr.7}E^i_t \ra E^\prime_t\eeq
 in $L^1_\text{loc}$, provided that
$\hnz\{U=t\}=0$. Thus (\ref{pr.7}) holds for all $t$ up to a countable
set $S=\{t\in [0,T]\ |\ \hnz\{U=t\}> 0\}$. Taking the limit we have
\beqs 
|\hnz(E^\prime_{t_1}\cap \Omega^\prime)-\hnz(E^\prime_{t_2}\cap
\Omega^\prime)| \leq C |t_2-t_1|^\frac{1}{k+1} 
\eeqs
for all $t_1,t_2 \in \R\setminus S$. Now let $t_0 \in S$ and pick
sequences $t_j^-\nearrow t_0,\ t_j^+\searrow t_0$ where $ t_j^-,t_j^+ \in
\R\setminus S$. Since 
$$ E^\prime_{t_j^-} \ra \{U\geq t_0\}\ , \ \ \ E^\prime_{t_j^+} \ra
\{U> t_0\} $$
this implies that $\hnz\{U=t_0\} =0$, and thus $S=\emptyset$. 
\epf 

We have seen before that the sets $E^i_t$ are minimizing area from outside
in $\Omega \times \R$. We now want to show that this property passes
to limit, even more we show that this is always preserved under
$L^1_\text{loc}$-convergence.

\bl\lb{pr.l5} Let $U \subset \rne$ be  open and $E_h\subset U$ a 
sequence of Caccioppoli-sets in $U$, which converge in
$L^1_\text{loc}(U)$ to $E\subset U$ such that $|\p^*E_h\cap K|\leq C(K)$ for all
$K\subset U$, $K$ compact, independently of $h$. If all the $E_h$ are minimizing
area from outside in $U$ then so does $E$.
\el
\bpf Let $E\subset F$ with $F\setminus E \subset K \subset U, K$
compact. Since $F\cup E_h \ra E$ and $E_h \ra E$ in
$L^1_\text{loc}(U\setminus K)$ there exists a compact set $K^\prime
\subset U$ with $K \subset \text{int}(K^\prime)$ with smooth boundary
$\p K^\prime$ such that
\beq\lb{pr.9}\begin{split}
\ \ |\p^*(F\cup E_h)\cap \p K^\prime|&\ =\ |\p^*(F\cap E_h)\cap \p
K^\prime| = |\p^*E_h \cap \p K^\prime| = 0\\  
\text{for all}\ h\ \ &\text{and}\ \    
\inl_{\p K^\prime}\! |\varphi^-_{F\cup E_h} -\varphi^+_{E_h}|\, d\hn \ra 0\ .
\end{split}\eeq 
Here $\varphi^-_{F\cup E_h},\ \varphi^+_{E_h}$ denote the inner,
resp. outer, trace of $F\cup E_h$ and $E_h$ on $\p K^\prime$. We can
assume w.l.o.g that $|\p^*E_h\cap U| \leq C$ for all $h$, and we
obtain for $F_h:=E_h\cup(F\cap K^\prime)$, compare \cite{Giusti},
Prop. 2.8:
$$
|\p^*F_h\cap U|=|\p^*E_h\cap(U\setminus K^\prime)|+\inl_{\p K^\prime}\!|\varphi^-_{F\cup E_h}
-\varphi^+_{E_h}|\, d\hn + |\p^*(F\cup E_h)\cap K^\prime|\ .
$$
The set $F_h$ is a valid comparison set for $E_h$, thus $|\p^*E_h\cap
U|\leq |\p^*F_h\cap U|$, which yields
$$|\p^*(F\cup E_h)\cap K^\prime| \geq |\p^* E_h \cap K^\prime| -
\inl_{\p K^\prime} \! |\varphi^-_{F\cup E_h}
-\varphi^+_{E_h}|\, d\hn\ .$$
Recall the general inequality
\beq\lb{pr.9a} |\p^*(E_1\cup E_2)\cap A| + |\p^*(E_1 \cap E_2)\cap A| \leq
|\p^*E_1\cap A| + |\p^*E_2 \cap A|\ , \eeq
which holds for any two Caccioppoli-sets $E_1,E_2$ in $U$, and $A$ any
open subset of $U$. By (\ref{pr.9}) we can apply this with $E_1=E_h,\
F=E_2$ and $A = K^\prime$ to get
$$|\p^*F\cap K^\prime| \geq |\p^* (E_h\cap F) \cap K^\prime| -
\int_{\p K^\prime}  |\varphi^-_{F\cup E_h}
-\varphi^+_{E_h}|\, d\hn\ .$$
Since $E_h\cap F \ra E$ in $L^1_\text{loc}$ we can use lower
semicontinuity and (\ref{pr.9}) to pass to limits:
$$|\p^*F\cap K^\prime| \geq |\p^*E\cap K^\prime|\ .$$
\epf  

\bc \lb{pr.c2} The sets $E^\prime_t$ are minimizing area from the
outside in $\Omega \times \R$ for all $t \in (0,T)$. As well the sets
$E_t$ are minimizing area from the outside in $\Omega$ for all $t \in
(0,T)$.
\ec
\bpf The proof of Lemma \ref{pr.l5} also works if we replace $\rne$ by
$N \times \R$. The first statement follows from Lemma \ref{pr.l4}. For the second statement let $F$ be a valid comparison set
for $E_t$ in $\Omega$, i.e. $E_t\subset F,\ F\setminus E_t \subset K \subset \Omega$, $K$
compact. Define $F^\prime:= (F\times (-l,l)) \cup
E^\prime_t$ which is a valid comparison set for $E^\prime_t$. Thus for
$K^\prime = K \times [-l+1,l+1]$ we have $|\p^*E^\prime_t\cap
K^\prime|\leq |\p^*F^\prime\cap K^\prime|$, i.e.
$$ 2l|\p^*E_t \cap K| \leq 2l|\p^*F\cap K| + 2\hne(F\setminus
E_t)\ .$$
Taking the limit $l\ra \infty$ proves the second statement.
\epf 
\bc\lb{pr.c3}
The function $t\mapsto |\p^* E_t|,\ t\in [0,T),$ is monotonically decreasing.
\ec

\bigskip

In the following we want to show the convergence of the Hypersurfaces
$$N^i_t \ra \Gamma_t\times \R$$
for a.e. $t$ in the sense of measures, where 
$$\Gamma_t:=\p\{u>t\} \subset \Omega \subset N , $$
i.e. $\Gamma_t\times \R = \p\{U>t\}$. So define Radon measures on
$\Omega \times \R$ by
$$ \mu^i_t:= \hne \res N^i_t\ , \ \ \ \ \mu_t:=
\hne\res\p^*E^\prime_t\ .$$
To prove the convergence $\mu^i_t \ra \mu_t$ we exploit the property
that the sets $E^i_t$ minimize area from the outside. We first want to
define a set $B\subset [0,T]$ of times where we can expect that such a
convergence holds true. Observe that we have $\p^*E_t \subset
\Gamma_t\subset \{u=t\}$ for all $t$. 

\bl\lb{pr.l6} There is a set $B\subset [0,T]$ of full measure, s.t. 
$$\hn(\{u=t\}\setminus \p^*E_t) = 0 $$
for all $t \in B$.
\el
\bpf Since $u$ is in $C^{0,1}(\bar{\Omega})\subset BV(\Omega)$ we can
compare the Coarea-formula for $BV$-functions and Lipschitz-functions
to get
$$\int_0^T\hn(\p^*E_t)\, dt = \int_\Omega |Du|\, d\hne = \int_0^T
\hn(\{u=t\})\, dt .$$
Since the integrals are finite, this yields
$$ \int_0^T \hn(\{u=t\}\setminus \p^*E_t)\, dt = 0\ ,$$
which implies the statement.
\epf

Thus $\mu_t = \hne \res (\Gamma_t\times \R)$ for all $t\in B$. Even
more $\mu_t$ is $(n+1)$-rectifiable for all $t \in B$. 

\bp\lb{pr.p1} For all $t \in B$, $\mu_t^i \ra \mu_t$ in the sense of
Radon measures.
\ep

\bpf We give the proof only in the case $N = \rne$. Since it
uses only local techniques it is straightforward to see that
the same
proof, with some minor modifications,  works also for a general $N$.\\ 
Fix a $t\in B$. By the mass bound (\ref{pr.4}) we can extract a
subsequence $\mu^{i_j}_t$ such that $\mu^{i_j}_t \ra \mu$, where $\mu$
is a Radon measure on $\Omega\times \R$.\\[+0.5ex]
\underline{Claim 1:}\ \ supp$(\mu) \subset \{u=t\}\times \R$.\\[+0.5ex]
Let $x \in \Omega \times \R,\ x \not\in \{u=t\}\times \R$,
i.e. $U(x)\neq t$. Let us assume $U(x) > t$. Thus there is a
$\delta>0$ such that $B_\delta(x)\Subset \Omega\times \R$ and $U(y)>t$
for all $t \in \bar{B}_\delta(x)$. So for $i$ sufficiently large
$U^{\eps_i} > t$ on $B_\delta(x)$, i.e. $N^i_t\cap B_\delta(x) =
\emptyset$ which implies $\mu^i_t(B_\delta(x))=0$ for large enough
$i$. So $x \not \in \text{supp}(\mu)$.\\[0.5ex]
\underline{Claim 2:}\  Let $B_\rho (x) \Subset \Omega\times \R$. Then
$$ \frac{\mu\big(\bar{B}_\rho(x)\big)}{\omega_{n+1}\rho^{n+1}} \leq
(n+2) \frac{\omega_{n+2}}{\omega_{n+1}}\ .$$
We have $\mu(B_\rho(x)) \leq \liminf_{j\ra \infty}
\mu^{i_j}_t(B_\rho(x))$. Using that  $\mu^{i_j}_t(B_\rho(x)) =$
\mbox{$|\p^*E^{i_j}_t\cap B_\rho(x)|$} and the $E^{i_j}_t$ minimize area from
the outside, we obtain by comparison with $E^{i_j}_t\cup B_\rho(x)$:
$$\mu^{i_j}_t(B_\rho(x)) \leq (n+2)\, \omega_{n+2}\rho^{n+1} $$
Thus for $\eps>0$:
$$\mu\big(\bar{B}_\rho (x)\big) \leq \mu(B_{\rho+\eps}(x)) \leq (n+2)\, 
\omega_{n+2}(\rho+\eps)^{n+1}\ ,$$
 which proves the claim for $\eps \ra 0$.\\[1.0ex]
The second claim establishes that $\mu$ is absolutely continuous w.r.t
to $\hne$-measure. By claim 1, Lemma \ref{pr.l6}, and the differentiation theorem for Radon
measures there is a function $\theta \in L^\infty(\Gamma_t\times \R\,
,\, \hne)$ such that we can write
\beq\lb{pr.10}\mu = \mu_t\res \theta \ .\eeq
\underline{Claim 3:}\ \  $\theta \geq 1$\ \ $\hne$-a.e. on
$\Gamma_t\times \R$.\\[+0.5ex]
By the differentiation theorem
\beq\lb{pr.11} \theta(x) = \lim_{\rho \ra 0}
\frac{\mu(B_\rho(x))}{\mu_t(B_\rho(x))}\ ,
\eeq
$\hne$-a.e. Let $x \in \Gamma_t\times \R$ such that this holds. Now
for all but at most countably many $\rho$, provided $B_\rho(x)\Subset
\Omega \times \R$, we have
$$ \mu(B_\rho(x)) = \lim_{j\ra \infty}\mu^{i_j}_t (B_\rho(x))\ .$$
On the other hand $\mu_t^{i_j} = \hne\res \p^*E^{i_j}_t$ and by lower semicontinuity
$$ \mu_t(B_\rho(x)) = |\p^*E^\prime_t\cap B_\rho(x)| \leq \liminf_{j
  \ra \infty}|\p^*E^{i_j}_t\cap B_\rho(x)|\ .$$ 
Thus $\theta \geq 1$ $\hne$-a.e. on $\Gamma_t\times \R$.\\[+0.5ex]
\underline{Claim 4:}\ \  $\theta(x) \leq 1$ for almost all $x \in
\p^*E^\prime_t$. \\[+0.5ex]
Here we use again the property that the sets $E^{i_j}_t$ are minimizing
area from the outside. Let $x \in \p^*E^\prime_t$ such that
(\ref{pr.11}) holds. By a translation we can assume that $x=0$. Since
$x \in \p^*E^\prime_t$ we know that as $\lambda \ra 0$,  the rescalings
$\lambda^{-1}\p^*E^\prime_t \ra T_x\p^*E^\prime_t$ in the sense of Radon measures,
and $\lambda^{-1}E_t^\prime \ra H$ in $L^1_\text{loc}$ where $H$ is
one of the halfspaces bounded by $T_x\p^*E^\prime_t$. By a rotation we
can assume that $T_x\p^*E^\prime_t = \{x^{n+2}=0\}$ and $H =
\{x^{n+2}< 0\}$. Let $\eps>0$ be given, and choose a $\lambda>0$
such that
\beq\lb{pr.12}
\bigg| \theta(x)-\frac{\mu(B_\lambda(x))}{\omega_{n+1}
  \lambda^{n+1}}\bigg|\leq \eps\ .
\eeq We can do this since
$\Theta^{n+1}(\p^*E^\prime_t, x) = 1$ and by (\ref{pr.11}). Then
define the $\delta$-slab 
$$S_\delta:= \{x \in \rnz\ | \ |x_{n+2}| \leq \delta\} \ .$$
By adjusting $\lambda$ maybe even further, we can assume that
\beq\lb{pr.13}
\hne(\lambda^{-1}\p^*E^\prime_t \cap (B_2\setminus S_\eps)) \leq
\eps \frac{\omega_{n+1}}{(n+2)\omega_{n+2}}
\eeq
and
\beq\lb{pr.14} 
\hnz\big((\lambda^{-1}E_t^\prime \setminus H)\cup (H \setminus
\lambda^{-1}E_t^\prime)\big)\leq \eps^2
\eeq
By claim 1 we have $\theta \leq (n+2)\omega_{n+2}/\omega_{n+1}$,
$\hne$-a.e. and thus (\ref{pr.13}) implies
\beq\lb{pr.15}
\mu_{\lambda}(B_2\setminus S_\eps) \leq \eps
\eeq
where $\mu_\lambda$ is the rescaling of $\mu$ by the factor
$\lambda^{-1}$,  defined by
$\mu_\lambda(A) = \lambda^{-(n+1)}\mu(\lambda A)$. We then choose $i$
(dropping the subscript $t$)
big enough such that 
\beq\lb{pr.16}
\big|\mu^{i_j}_\lambda(B_2\setminus S_\eps) - \mu_\lambda(B_2\setminus
S_\eps)\big| \leq \eps\ , \ \ \ \ \ \ \big|\mu^{i_j}_\lambda(
B_2) - \mu_\lambda( B_2)\big| \leq \eps
\eeq
and 
\beq\lb{pr.17}
\hnz\big((\lambda^{-1}E_t^\prime \setminus \lambda^{-1}E_t^{i_j})\cup
(\lambda^{-1}E_t^{i_j}\setminus \lambda^{-1}E_t^\prime)\big) \leq
\eps^2
\eeq
Combining this with (\ref{pr.15}) and (\ref{pr.14}) we get
\beq\lb{pr.18}
\mu^{i_j}_{\lambda}(B_2\setminus S_\eps) \leq 2\eps
\eeq
and
\beq\lb{pr.19}
\hnz\big((\lambda^{-1}E_t^{i_j} \setminus H)\cup (H \setminus
\lambda^{-1}E_t^{i_j})\big)\leq 2 \eps^2\ .
\eeq
In other words, up a set of small measure, $\lambda^{-1}E^{i_j}_t$
looks like the halfspace $H$ on $B_2$. We now employ that
$\lambda^{-1}E^{i_j}_t$ is minimizing area from the outside to get an
upper bound on the area contained in the slab $S_\eps\cap B_1$. Now by
(\ref{pr.19}) there is a $\delta \in [\eps, 2\eps]$ such that
\beq\lb{pr.20}
\hne\Big(\big(\{x_{n+2}= - \delta\}\cap (\lambda^{-1}E^{i_j}_t)^C\big)\cap
B_2\Big) \leq 2\eps
\eeq
We take as comparison set $F$,
$$ F:= \lambda^{-1}E^{i_j}_t \cup \big(S_\delta\cap B_1)\ .$$
Since $\lambda^{-1}E^{i_j}_t$ is minimizing area from the outside,
together with (\ref{pr.20}), (\ref{pr.18}), this implies
\beqs
\mu^{i_j}_\lambda (B_1) \leq \omega_{n+1}+ 4 \eps n \omega_{n+1} + 4
\eps \ ,
\eeqs
which in turn gives with (\ref{pr.16}) that
$$ \mu_\lambda(B_1) \leq \omega_{n+1}+ 4 \eps n \omega_{n+1} + 5\eps\ $$
Finally applying this to (\ref{pr.12}), we arrive at
$$ \theta(x) \leq 1 + (4n + 5/\omega_{n+1} + 1)\, \eps \ .$$
which proves the claim.\\[+1ex] 
To finally prove that $\mu = \mu_t$ we combine claim 3 and claim 4 to
see that $\theta = 1 $ $\hne$-a.e. on $\Gamma_t \times \R$, and use
(\ref{pr.10}). Thus the limit measure $\mu$ does not depend on the
subsequence, so the whole sequence converges, i.e. $\mu^i_t \ra \mu_t$. 
\epf

Recall that we have the uniform space-time bound, see (\ref{pr.5}), 
$$ \inl_{0}^T\!\! \int_{\Omega^\prime} H_i^{k+1}\, d\mu^i_t\, dt \leq
(b-a)\hn(\p\Omega)+ 2 \hne(\Omega)\ ,$$
where $\Omega^\prime:=\Omega\times (a,b)$. Thus by Fatou's lemma
\beq\lb{pr.22}\liminf_{i\ra \infty}\int_{\Omega^\prime} H_i^{k+1}\, d\mu^i_t <
\infty\eeq
for almost all $t \in [0,T]$. Then we let
$$ \tilde{B}:= \{t\in B\ |\ \text{(\ref{pr.22}) holds for any 
bounded interval}\ (a,b)\}. $$
Note that $\tilde{B}$ again has full measure. Thus for every $t\in
\tilde{B}$ there is subsequence $(i_j)$ such that by the mass bound
and H\"older's inequality
\beq \lb{pr.23} \sup_{j\geq 0} \int_{A} |H^{i_j}|\, d\mu^{i_j}_t \leq C
  \big(\mu^{i_j}_t(A)\big)^\frac{k}{k+1} < C(\Omega) 
\eeq
for every $A \subset \Omega^\prime$. By the compactness theorem
for $(n+1)$-rectifiable varifolds of Allard, there is a further subsequence
(which we again denote by $(i_j)$) such that $N^{i_j}_t$ converges in the sense
of varifolds to a limit, which again is $(n+1)$-rectifiable. Since
$\mu^i_t \ra \mu_t$ this implies
$$ N^{i_j}_t \ra \Gamma_t\times \R $$
in the sense of varifolds, where we see $\Gamma_t\times \R$ as a
$(n+1)$-rectifiable, unit density varifold. The estimate (\ref{pr.23})
then implies that the total variation $\delta (\Gamma_t\times \R)$ is
absolutely continuous w.r.t. $\mu_t$, i.e. $\Gamma_t\times \R$ carries
a weak mean curvature $\mathbf{H}$, and thus by the product structure
also $\Gamma_t$. The varifold convergence now
implies that
$$ \mu^i_t \res (-H^i\nu^i_t) \ra \mu_t \res \mathbf{H} $$
in the sense of vector valued Radon measures. So by lower
semicontinuity-results  for convex functionals of
Hutchinson, see \cite{Hutchinson},
$$ \int_{\Omega^\prime} |\mathbf{H}|^{k+1}\, d\mu_t \leq \liminf_{i
  \ra \infty} \int_{\Omega^\prime} H_i^{k+1}\, d\mu^i_t\ .$$
So again by Fatou's lemma
\beqs\begin{split}\inl_{0}^T\!\!\! \int_{\Omega\times(a,b)} |\mathbf{H}|^{k+1}\, d\mu_t\, dt
\leq&\ \liminf_{i \ra \infty}\inl_{0}^T\!\!\!\int_{\Omega^\prime}
H_i^{k+1}\, d\mu^i_t\, dt\\
 \leq&\ (b-a)\hn(\p\Omega)+ 2 \hne(\Omega)\
. 
\end{split}
\eeqs
This implies
\beq\lb{pr.24}
\inl_0^T\!\!\!\int_{\Omega\cap \Gamma_t} |\mathbf{H}|^{k+1}\, d\hn\, dt \leq
\hn(\p\Omega)\ .
\eeq
Applying Allard's regularity theorem we can summarize the above in the following
regularity and approximation result.

\bt\lb{pr.t1} There is a set $\tilde{B}\subset [0,T]$ of full measure,
such that for all $t\in \tilde{B}$ the following is true.
\begin{itemize}
\item[i)] \ \ $\hn(\Gamma_t\setminus \p^* E_t) = 0$
\item[ii)] For $k>n-1$ the surfaces $\Gamma_t$ are up to a closet set
  $A_t\subset \Gamma_t$ with $\hn(A_t)=0$ in $C^{1,1-\frac{k+1}{n}}$.
\item[iii)] There is a subsequence $(i_j)$, depending on $t$, such that
  $N^{i_j}_t \ra \Gamma_t\times \R$ in the sense of varifolds. If
  $k>n$, then away
  from the set $A_t\times \R$, this convergence is in $C^{1,\alpha}$
  for any $0<\alpha<1-\frac{k+1}{n+1}$.
\end{itemize}
\et

\medskip

\section{The main estimate}\lb{mainestimate}

In this section we show that the estimate (\ref{intro.2}) is
valid for almost all $t \in [0,T]$ if $N=\R^{n+1}$. The estimate in
case $N$ is a complete, simply-connected $3$-manifold with nonpositve
sectional curvatures will be given at the end of the section.

\medskip

Our aim is to transfer the computation presented in (\ref{sm.2})
for the smooth case to the setting of lower regularity of
solutions of the weak flow. Assume that $k>n$. By Theorem \ref{pr.t1}
and Lemma \ref{pr.l3}
we know that there is a set $\tilde{B} \subset [0,T]$ of full measure, such
that for $t\in \tilde{B}$ the following statements are true: Up to a closed set of
$\hn$-measure zero, the set $\Gamma_t=\p\{u>t\}$ is a
$C^{1,\alpha}$ hypersurface, which
carries a weak mean curvature in $L^{k+1}(\Gamma_t)$. We can
as well assume that there is a sequence $\eps_i \ra 0$ such that
\beqs
N_t^{i} \ra \Gamma_t \times \R
\eeqs
in the sense of varifolds, and
\beq\lb{me.2}
\limsup_{i\ra \infty} \inl_{\neit \cap (\Omega\times I)} \!\!\!\!\!\! |H_i|^{k+1}+ |\nabla
H_i|^2 \, d\meit\  < \ \infty\ ,
\eeq
for $I$ a bounded interval. Our aim is to do all the
computations on equidistant hypersurfaces to $\Gamma_t$ and then
pass to limits. For this purpose, define for $s>0$
$$
L_s(\Gamma_t) = \p\{x\in \rn \ |\  \text{dist}(x,E_t)\leq s\}\ .
$$

\bl\lb{me.l1}Let $\phi \in C^0_c(\R^{n+2})$. Then for all $p<k+1$
\beq\lb{me.2b}
\int_{\neit} \phi H_i^p\, d\hne\ \ra \int_{\Gamma_t\times \R}\phi |\mathbf{H}|^p
\, d\hne\ .
\eeq
\el
\bpf By the results stated above we know that 
$$\neit \ra \Gamma_t \times \R$$
not only in the sense of varifolds, but also by Allard's theorem away
from the singular set sing$(\Gamma_t)\times \R$ locally uniformly in
$C^{1,\alpha}$. Furthermore sing$(\Gamma_t)\times \R$ is closed and has
$\hne$-measure zero. So given any $\delta >0$ there is a neighborhood
$S$ of sing$(\Gamma_t)\times \R$ such that
\beq\lb{me.2c}
\limsup_{i\ra \infty}\bigg|\int_{\neit \cap S}\phi H^p\, d\hne \bigg|
\leq C \sup |\phi|\  \delta\ .
\eeq
Here we used the uniform $L^{k+1}$-estimate on $H$ from
(\ref{me.2}). Now outside $S$ the convergence of the hypersurfaces is
locally uniform in $C^{1,\alpha}$. So it suffices to check
(\ref{me.2b}) locally, such that $\Gamma_t\times \R$ and $\neit$ can
be written as converging graphs (with bounded gradient) over a fixed
hyperplane. Since the hypersurfaces converge as varifolds, and we have
local convergence in $C^1$, the mean curvature $H_i$ converges weakly to $H$. The
uniform $L^2$-estimate on $\nabla H_i$ from (\ref{me.2}) gives
that $H_i \ra H$ in $L^2$, and by interpolation in $L^p$ for every
$p<k+1$. Using a suitable partition of unity and (\ref{me.2c}) we
get the claimed convergence.
\epf

To be able later to control the convergence as $s\ra 0$ we need to
control the local area growth of the hypersurfaces $L_s$ in the
parameter $s$.

\bl\lb{me.l2} For almost all $s>0$ the following statement is true:
Let $K\subset L_s(\Gamma_t)$ be \hn-measurable and define
$P(K):=\{x\in \Gamma_t\ |\ \exists\ y \in K\ \text{with}\
|y-x|=s\}$. Then

$$ \hn(K) \leq \frac{1}{n^n}\inl_{P(K)\cap \Gamma_t} \!\!\!\!\!\!
(|\mathbf{H}(x)|\, s+n)^n\, d\hn(x)\ . $$
\el

\bpf We can assume that $K$ is compact, and so $P(K)$ is compact as
well. Define for $\gamma,\eta>0$
$$ K_\gamma:= \{x \in \rne\ |\ \text{dist}(x,K)< \gamma\}\ , \ \ P_\eta:=
\{x \in \rne\ |\ \text{dist}(x,P(K))< \eta \}\ .
$$

\medskip

\noindent\underline{Claim 1:} For all $\eta > 0$ there exists a $\gamma_0 >0$ such
that $P(K_\gamma\cap L_s(\Gamma_t)) \subset P_\eta$ for all $\gamma\leq\gamma_0$. \\[+1.0ex]
Assume to the contrary that there are points $y_i\in
L_s(\Gamma_t)\setminus K $ with $y_i \ra y_\infty \in K$ and points
$x_i \in \Gamma_t$ with $|y_i-x_i|=s$, but dist$(x_i, P(K))\geq
\eta$. We can assume that $x_i \ra x_\infty \in \Gamma_t$. Thus
$|y_\infty-x_\infty|=s$, but $x_\infty \not \in P(K)$. This proves
Claim 1.\\
\medskip
Since the distance function to $E_t$ is lipschitz, we can argue as in
Lemma \ref{pr.l6}
that for a.e. s the set $\{x \in \Omega\ | \ \text{dist}(x,E_t) < s
\}$ is a Caccioppoli-set and
\beq\lb{me.3} 
\p^*\{x \in \Omega\ | \ \text{dist}(x,E_t) < s
\}= L_s(\Gamma_t) 
\eeq
up to $\hn$-measure zero. A further thing to note is that by
(\ref{me.2}) the mean curvature $H^{\eps_i}$ is uniformly bounded in
$L^p(\neit)$ for some $p>n+1$. This gives that also 
\beqs
\neit \ra \Gamma_t\times \R
\eeqs
locally in Hausdorff-distance, which in turn implies that
$$\{z \in \Omega\times \R | \ \text{dist}(z, \tilde{E}^{i}_t) < s
\}\ra \{x \in \Omega\ | \ \text{dist}(x,E_t) < s
\}\times \R$$
in $L^1_\text{loc}$. Then using the lower semicontinuity of the
BV-norm we deduce, note (\ref{me.3}),
\beq\lb{me.5}
\hne((K_\gamma\times I)\cap (L_s(\Gamma_t\times \R))) \leq \liminf_{i\ra \infty}
\hne((K_\gamma \times I) \cap L_s(\neit))\ ,
\eeq
for any bounded interval $I$. We now want to apply  a result of Li and Nirenberg
\cite{LiNirenberg05} or equivalently Itoh and Tanaka
\cite{ItohTanaka01}, which says: Given a bounded open set $S\subset
\R^{n+2}$ with smooth boundary, define $G$ to be the largest open subset of
$S$ such that every point $x$ in $G$ has a unique closest point on $\p
S$. Then the set $\Sigma(S):=
S\setminus G$ has finite $\hne$-measure. Furthermore for every $x \in
G$ the distance function to the boundary is smooth. \\
\medskip
Since the sets $\R^{n+2}\setminus\tilde{E}^{i}_t$ have a smooth
boundary and converge locally in Haussdorff-distance to
$\R^{n+2}\setminus (E_t\times \R)$ we can apply the above result to
deduce that the sets $\Sigma_i\subset
\R^{n+2}\setminus\tilde{E}^{i}_t$, defined as above have locally
finite $\hne$-measure. Thus for almost all $s\in (0,1)$ 
\beq\lb{me.6}
\hne(L_s(\neit)\cap \Sigma_i) = 0\ \ \text{for all}\ i\in \N\ .
\eeq
Now pick an $s_0\in (0,1)$ such that (\ref{me.5}) and (\ref{me.6})
hold. Let $x_0\in L_{s_0}(\neit)\setminus \Sigma_i$. Then there is
neighborhood of $x_0$ such that every point has a unique closest
point on $\neit$ and the distance to $\neit$ is smooth. Even more, the
same is true for a neighborhood of the line $l$ connecting $x_0$ to its
closest point on $\neit$. We compute along $l$
\beq\lb{me.7} 
\frac{\p}{\p s} H = -|A|^2 \leq -\frac{1}{n}H^2
\eeq
where $H$ is the mean curvature of $L_s(\neit)$ along $l$. Comparing
with the solution of the ODE yields
\beq\lb{me.8}
H(s)\leq \max\Big(\frac{n H(0)}{H(0)\,s + n}, 0 \Big) = \frac{n
  H(0)}{H(0)\,s + n}\ ,
\eeq
where the last equality holds since $H(0)>0$. The evolution of the
measure along $l$ is given by
\beqs 
\frac{\p}{\p s} d\mu = H(s)\, d\mu
\eeqs
which can be integrated to
\beq\lb{me.10}
d\mu(s)= \exp\bigg(\int_0^s H(\tau)\, d\tau\bigg)\, d\mu(0)\ .
\eeq
Inserting the estimate (\ref{me.8})
\beqs
\int_0^s H(\tau)\, d\tau \leq \int_0^s\frac{n H(0)}{H(0)\, \tau +n}\,
d\tau = \log\bigg( \bigg(\frac{H(0)\,s+n}{n}\bigg)^n\bigg)\ ,
\eeqs
we arrive at
\beq\lb{me.12}
d\mu(s)\leq \bigg(\frac{H(0)\,s+n}{n}\bigg)^n  d\mu(0)\ .
\eeq 
If we denote with $W_\gamma^i=\{z\in\neit\ |\ \exists\  w \in 
(K_\gamma \times (0,1)) \cap L_s(\neit)\ \text{with}\ |z-w|=s \}$ we
estimate, using (\ref{me.6}) and (\ref{me.12}),
\beq\lb{me.13}
\hne((K_\gamma \times I) \cap L_s(\neit))\leq
\frac{1}{n^n}\int_{W_\gamma^i \cap \neit} (H s + n)^n\, d\hne
\eeq
for all $i \in \N$. Given an $\eta > 0$ assume that $\gamma>0$ is small
enough such that also $2\gamma < \gamma_0$. \\[+2ex]
\underline{Claim 2:} For $i$ large enough
$W_\gamma^i \subset P_{2\eta}\times (-\eta, 1+\eta)$.\\[+2ex]
Assume this would not be the case. Then there exist points $w_i \in 
(K_\gamma \times (0,1)) \cap L_s(\neit)$ and points $z_i \in \neit$
with $|w_i-z_i|=s$ and dist$(z_i, P_\eta\times(0,1))\geq \eta$. We can
further assume that $w_i \ra w_\infty$ and $z_i \ra z_\infty$ with
$|w_\infty-z_\infty|=s$. 
Since $\neit \ra \Gamma_t\times \R$ in Hausdorff-distance we have that
$$w_\infty \in (K_{2\gamma} \times [0,1])\cap (L_s(\Gamma_t)\times
[0,1])\ \ 
\text{and}\ \ z_\infty\in\Gamma_t$$ 
but
$$ z_\infty \not \in (P_{\eta}\times [0, 1])\cap (\Gamma_t\times \R)\ .$$
This contradicts Claim 1.\\
\medskip
We now combine (\ref{me.5}), (\ref{me.13}), Claim 2, and Lemma \ref{me.l1} to arrive
at
\beqs
\hne((K_\gamma\times (0,1))\cap (L_s(\Gamma_t\times \R))) \leq 
\inl_{(P_{2\eta}\times (-\eta, 1+\eta)) \cap
  (\Gamma_t\times\R)} \!\!\!\!\!\!\!\!\!\!\!\!\! \!\!\!\!\!\!\!\! n^{-n}(H s + n)^n\, d\hne\ ,
\eeqs
for almost all $\eta$, and $\gamma$ chosen appropriately. Then let
$\eta,\gamma \ra 0$.
\epf

The next lemma shows that there can't be two different points $p_1,p_2
\in L_s(\Gamma_t)$ and a point $x_0 \in \Gamma_t$ with $|p_1-x_0|=|p_s-x_0|=s$.  

\bl\lb{me.l3} Assume that in a point $x_0\in \Gamma_t$ the set $E_t$
can be touched from outside by two balls $B_s(p_1), B_s(p_2)$ for $s>0$. Then $p_1=p_2$.
\el
\bpf
Assume to the contrary that two different balls $B_s(p_1), B_s(p_2)$ touch
$E_t$ in $x_0$. Since $\mathbf{H} \in L^p(\Gamma_t)$
for some $p>n$ the density $\Theta^n$ exists at $x_0$. Since $E_t$
minimizes area from the outside, we can argue as in the proof of
Proposition \ref{pr.p1} to show that $\Theta^n$ satisfies the
bound $\Theta^n(x_0) \leq C(n)$ . By upper semicontinuity
$\Theta^n(x_0)\geq 1$.  Thus the blow-ups
$$ \frac{1}{\lambda_i}(\Gamma_t-x_0) $$
converge for some sequence $\lambda_i\ra 0$ to a stationary cone $C$
with $\Theta^n_C(0)=\Theta^n_{\Gamma_t}(x_0)$. Since $p_1\neq p_2$ the
cone $C$ has to be a subset of two different closed halfspaces
$T_1,T_2$ where $0$ is contained in either of the boundaries of
$T_1,T_2$. This implies that supp$(C)$ has to be a subset of $\p T_1$ (see for
example \cite{Simon}, Thm. 36.5). Since the same holds also for $\p T_2$
we have $C \subset \p T_1 \cap \p T_2$. This yields $\hn(C)=0$, which
contradicts $\Theta_C^n(0)>0$.
\epf

Given a $\hn$-measurable function $f$ on $L_s(\Gamma_t)$ this enables us to
define the ``pull-back''  $\tilde{f}$ to $\Gamma_t$ by
\beqs
\tilde{f}(x) \ =\ \begin{cases}
\ f(y)\ \text{if}\ \exists \ y\in L_s(\Gamma_t)\ \text{with}\ |y-x|=s,\\
\ 0\ \ \text{else},
\end{cases}
\eeqs
for any $x \in \Gamma_t$. By the regularity of $\Gamma_t$, $\tilde{f}$
is also $\hn$-measurable and we get the following Corollary:

\bc\lb{me.c4}
For almost all $s>0$ the following statement is true:
Let $f$ be a $\hn$-measurable, nonnegative function on $L_s(\Gamma_t)$ and define
$\tilde{f}$ as above. Then

$$ \inl_{L_s(\Gamma_t)}f(y)\, d\hn(y) \leq \frac{1}{n^n}\inl_{\Gamma_t} 
\tilde{f}(x)\ (|\mathbf{H}(x)|\, s+n)^n\, d\hn(x)\ . $$
\ec

In the next proposition we use the previous estimate to prove
(\ref{intro.2}) for almost every $t$.

\bp\lb{me.p5} Let $k>n$. Then for all $t\in \tilde{B}$ the estimate
\beq\lb{me.16}
\int_{\Gamma_t}|\mathbf{H}|^n\, d\mu \geq \Big(\frac{n}{n+1}\, c_n\Big)^n\ 
\eeq
holds.
\ep
\bpf We first need to investigate further the regularity of the
hypersurfaces $\Gamma_t$.\\[2ex]
\underline{Claim 1:} Let $t \in \tilde{B}$. Then $\Gamma_t$ is twice
differentiable $\hn$-a.e.\\[2ex]
Let $S\subset \Gamma_t$ be the singular part of $\Gamma_t$, such that
away from  $S$, $\Gamma_t$ can be written locally as the graph of a
$C^{1,\alpha}$-function $u$. By (\ref{me.2}) we know that the mean
curvature $H$ of $\Gamma_t$ is in $L^p(\Gamma_t)$ for some $p>n$. This
implies that $u$ is a weak solution of the equation 
$$\text{div}\bigg(\frac{Du}{\sqrt{1+|Du|^2}}\bigg) = H .$$ 
Since we can assume that $p\geq 2$ we can deduce as for example in
\cite{Giusti}, Theorem C.1, that
$u\in W^{2,2}$ and $u$ is a strong solution of the equation
\beq\lb{me.17}
a_{ij}(Du)\, D_{ij}u = \sqrt{1+|Du|^2} H\ ,
\eeq where
$$ a_{ij}(Du)=\delta_{ij}-\frac{D_iu\, D_ju}{1+|Du|^2}\ .$$
Since $u$ is in $C^1$ the coefficients $a_{ij}$ are uniformly elliptic
and continuous. Thus we can apply the estimates of Calderon-Zygmund to
deduce that $u\in W^{2,p}$. Since $p>n, u$ is twice differentiable
$\hn$-a.e. (see for example \cite{EvansGariepy}, section 6.4, Theorem
1). Note that  at a point $x_0$ where $u$ is twice differentiable, we can write
\beqs\begin{split}
u(x)=&\ u(x_0) + D_iu(x_0)(x-x_0)_i +\frac{1}{2} D_{ij}u(x_0)
(x-x_0)_i (x-x_0)_j\\
&\ \ \  + o(|x-x_0|^2)\ ,
\end{split}\eeqs 
and equation (\ref{me.17}) holds pointwise a.e.\\[2ex]
We now argue in a similar spirit as \cite{Kleiner92}. Let us denote with
$C(\Gamma_t)$ the outer convex hull of
$\Gamma_t$ and let $C_s=L_s(C(\Gamma_t))$. Observe that forming the
convex hull and taking an outer equidistant surface commutes, i.e.
$$C_s=L_s(C(\Gamma_t))=C(L_s(\Gamma_t))\ .$$
$C_s$ is convex and $C^{1,1}$ (see Appendix B in
\cite{Almgren86}). The nearest point projection \mbox{$\pi_s:C_s\ra C_0$} is
well defined, distance nonincreasing and is a bijection between $C_s\cap
L_s(\Gamma_t)$ and $\Gamma_t\cap C_0$. \\[2ex]
\underline{Claim 2:} Let $s>0$. For $\hn$-a.e. $p \in \Gamma_t\cap C_0$ we can estimate
\beqs
\frac{n H_{\Gamma_t}\!(p)}{H_{\Gamma_t}\!(p)\,s + n}\geq
H_{C_s}(\pi_s^{-1}(p))\ .
\eeqs
We can assume that $p$ is not in the singular part of $\Gamma_t$ and,
by Claim 1, that $\Gamma_t$ is twice differentiable at $p$. Since
$\hn$-zero sets on $C_s$ are mapped under $\pi_s$ to $\hn$-zero sets
on $\Gamma_t\cup C_0$, we can as well assume that $C_s$ is twice
differentiable at $q:=\pi^{-1}_s(p)$. Now let $\Sigma$ be a supporting
hypersurface of $C_s$ at $q$, i.e. $\Sigma$ is locally a smooth
hypersurface, which touches $C_s$ in $q$ from the outside. For a given
$\delta >0$ we can assume that
\beq\lb{me.20}
H_{C_s}(q)-\delta\leq H_\Sigma(q)\leq H_{C_s}(q) .
\eeq
Furthermore, we may assume that $\Sigma$ is a convex hypersurface,
with principal curvatures local bounded around $q$ by $1/s$. We 
translate $\Sigma$ such that it touches $\Gamma_t$ in $p$ from the
outside. By Claim 1, this implies that
\beqs
H_{\Gamma_t}(p)\geq H_\Sigma(q)\ .
\eeqs
We now evolve $\Sigma$ (translated back to its original position)
equidistantly towards $\Gamma_t$. Let us denote with $\Sigma_\tau$ the
so obtained hypersurface with dist$(\Sigma,\Sigma_\tau)=\tau$ for
$0<\tau <s$. Since we have assumed that the principal curvatures of
$\Sigma $ are locally bounded by $1/s$, $\Sigma_\tau$ remains smooth,
and it touches $C_{s-\tau}$ in $q^\prime$ from the outside, with
$\pi_{s-\tau}(q^\prime)=p$. Now, as above, this yields
\beq\lb{me.22}
H_{\Gamma_t}(p)\geq H_{\Sigma_{\tau}}(q^\prime)\ .
\eeq
We use the evolution of the mean curvature under the equidistant flow,
 as in the proof of Lemma \ref{me.l2}, see (\ref{me.7}) and
 (\ref{me.8}), to obtain
\beq\lb{me.23}
H_\Sigma(q)\leq \frac{n
  H_{\Sigma_\tau}(q^\prime)}{H_{\Sigma_\tau}(q^\prime) \tau + n}\ .
\eeq
Now$$\frac{\partial}{\partial H} \bigg(\frac{nH}{H\tau+n}\bigg) \geq 0\
.$$
So we can combine (\ref{me.20}), (\ref{me.22}) and (\ref{me.23}) to
conclude that
\beqs
H_{C_s}(q)-\delta \leq \frac{n
  H_{\Gamma_t}(p)}{H_{\Gamma_t}(p) \tau + n}\ .
\eeqs
Taking the limits $\tau\ra s$ and $\delta \ra 0$ proves Claim
2.\\[2ex]
We finally prove the proposition. Since the
Weingarten map on $C_s$, restricted to $C_s\cap L_s(\Gamma_t)$, covers
$\mathbb{S}^n$ at least once we can estimate as in the smooth case
\beqs
|\mathbb{S}^n| \leq \inl_{C_s\cap L_s(\Gamma_t)}\!\!\!\!\!\!G \, d\hn \leq \frac{1}{n^n}
\inl_{C_s\cap L_s(\Gamma_t)}\!\!\!\!\!\!\big(H_{C_s}\big)^n\, d\hn \ .
\eeqs
We estimate further, applying Corollary \ref{me.c4} and Claim 2, to
arrive at
\beqs
\begin{split}
n^n|\mathbb{S}^n| \leq&\  \frac{1}{n^n}\inl_{\Gamma_t\cap
  C_0}\!\!\!\!\big(H_{C_s}(\pi_s^{-1}(x))\big)^n \big(H_{\Gamma_t}(x) s + n
\big)^n\, d\hn(x)\\
\leq& \inl_{\Gamma_t\cap C_0}\big(H_{\Gamma_t}\big)^n(x)\, d\hn(x)
\leq \inl_{\Gamma_t}\big|\mathbf{H}_{\Gamma_t}\big|^n\, d\hn\ .
\end{split}
\eeqs
\epf

The estimate (\ref{intro.2}) in the case that the ambient manifold is
a 3-dimensional Hadamard manifold $(N^3,\bar{g})$, is a 
modification of an argument due to L. Simon,
see \cite{Simon93}. The main tool is, as in the monotonicity formula,
to use a well chosen vectorfield in the the first variation identity
\beq\lb{me.27}
\inl_{M}\text{div}_{M}(Y)\, d\mu = -\inl_{M} \bar{g}(Y,{\mathbf H})\, d\mu
\eeq
which holds for any $C^{0,1}$ vector field $Y$, defined in a neighborhood of $M$. As
before we write for $p\in M$:
$$ \text{div}_M (Y)(p) = \sum_{i=1}^2
\bar{g}(\bar{\nabla}_{e_i}Y,e_i)\ ,$$
where $\bar{\nabla}$ is the covariant derivative on $N$ and
$e_1,e_2$ form an orthonormal basis of $T_pM$. Since we still can make
sense of (\ref{me.27}) in the varifold setting, the estimate needs
much less regularity than in case we treated before.

\bl\lb{me.l6} Let $N^3$ be a complete, simply connected
3-manifold, with nonpositive sectional curvatures. Let $M\subset
N^3$ be a bounded integer 2-rectifiable varifold, carrying a weak mean curvature
$\mathbf H \in L^2(\mu)$. Then
\beq\lb{me.28}
 \inl_M |\mathbf{H}|^2 \, d\mu \geq 16\pi \ . 
\eeq
\el   
\bpf For the computation in the case of a flat Euclidean ambient space
one uses the position vectorfield $X$, centered at a point $x_0$. The
calculation then
depends on the fact that
\beqs
\text{div}_{M}(X)(x) = 2\ ,
\eeqs
for all $x\in M,\ x \neq x_0$, such that the tangent space of $M$
exists at $x$. Now on a complete, simply-connected manifold
$N$, we replace this vectorfield by
\beqs
X:= r \bar{\nabla} r\ ,
\eeqs
where $r(p):=\text{dist}_{N}(p,p_0)$ for a fixed
$p_0\in N$. Here $\bar{\nabla}$ denotes the gradient
operator on $N$. Let us assume that the sectional
curvatures of $N$ are bounded above by $-\kappa$ for some
$\kappa\geq 0$. The distance function to a point on such a manifold has two important
properties, see for example \cite{Petersen}:
\beqs
\begin{split}
\bar{\nabla}r \neq&\ 0 \\
\text{Hess}(r)=\bar{\nabla}^2r\geq&\
\Psi(r)\big(\text{id}-\bar{\nabla}r\otimes\bar{\nabla}r\big) ,
\end{split}
\eeqs
for $p\neq p_0$, where
$\Psi(r):=\sqrt{\kappa}\cosh(\sqrt{\kappa}r)/\sinh(\sqrt{\kappa}r)$
and the second inequality holds w.r.t.\! an orthonormal basis of $T_pN$. For a
point $p\in M$, such that the tangent space of $M$ exists at $p$, choose a normal vector $\nu$ to $T_pM$ and compute
\beq\lb{me.32}
\begin{split}
\text{div}_{M}(X) =&\ \text{div}_M (r\bar{\nabla}r) =
r\,\text{div}_M(\bar{\nabla}r) + \bar{g}(\nabla_M r, \bar{\nabla}r)\\
=&\  r\, \text{tr}_{T_pM}\big(\text{Hess}(r)\big) + 1 -
\bar{g}(\bar{\nabla} r,\nu)^2\\
\geq&\ r\,\Psi(r)\,
\text{tr}_{T_pM}\big(\text{id}-\bar{\nabla}r\otimes\bar{\nabla}r\big)
+1 - \bar{g}(\bar{\nabla} r,\nu)^2\\
=&\ \Big(\text{tr}_{T_pN}\big(\text{id}-
\bar{\nabla}r\otimes\bar{\nabla}r\big)+\bar{g}(\bar{\nabla}
r,\nu)^2-1\Big)+1 - \bar{g}(\bar{\nabla} r,\nu)^2\\
&\ +\big(r\,\Psi(r)-1\big)\text{tr}_{T_pM}\big(\text{id}-\bar{\nabla}r\otimes\bar{\nabla}r\big)\\
=&\ 2+\big(r\,\Psi(r)-1\big)\Big(1+ \bar{g}(\bar{\nabla} r,\nu)^2\Big)
\geq 2\ ,
\end{split}
\eeq
since $r\Psi(r)\geq 1$. The computation now proceeds as in the Euclidean case: Pick any
$p_0\in M$ such that the density $\Theta(p_0)$ exists and
$\Theta(p_0)\geq 1$. For $0<\sigma<\rho$ we can substitute in
(\ref{me.27}) the vectorfield $Y(p)\equiv
(|X|_\sigma^{-2}-\rho^{-2})_+ X$ where $|X|_\sigma=\max (|X|,
\sigma)$. A direct computation, using (\ref{me.32}), then yields
\beqs
\begin{split}
2\sigma^{-2}\mu(B_\sigma)+ 2\!\!\!\!\!\! \inl_{B_\rho\setminus
  B_\sigma} \!\!\!\!\bigg|\frac{1}{4}\mathbf{H}+&\ \frac{X^\perp}{|X|^2}\bigg|^2
d\mu \leq  2\rho^{-2}\mu(B_\rho)+\frac{1}{8}\!\!\!\!\inl_{B_\rho\setminus
  B_\sigma}\!\!\!\!\!\!|\mathbf{H}|^2\, d\mu\\
&\ -\sigma^{-2}\inl_{B_\sigma}\bar{g}(X,\mathbf{H})\, d\mu +
\rho^{-2}\inl_{B_\rho}\bar{g}(X,\mathbf{H})\, d\mu\ ,
\end{split}
\eeqs  
where we assume that all balls $B_\sigma,B_\rho$ are centered at
$p_0$. Since 
$$\lim_{\sigma\ra 0}\sigma^{-2}\mu(B_\sigma) \geq \pi\ ,$$
 we can take the limits $\sigma \ra 0$ and $\rho \ra \infty$ to prove the
estimate (\ref{me.28}).\epf

In the case that $N^3$ has sectional curvatures bounded
above by $-\kappa$ for some $\kappa>0$ there is a stronger estimate
by the Gauss-Bonnet formula. This estimate was used by B. Kleiner in
\cite{Kleiner92} to show that on such a manifold $N^3$ the
isoperimetric inequality of the model space with constant sectional
curvatures $\kappa$ holds. The estimate goes as follows. Let $M\subset
N^3$ be a closed hypersurface, diffeomorphic to a sphere and
denote with $R_\text{int}$ the intrinsic scalar curvature of $M$. Then
by the Gauss equations
\beq\lb{me.34}
4\pi = \int_M R_\text{int}\, d\hz \leq \int_M G -\kappa\, d\hz \leq
\int_M \frac{1}{4} H^2 - \kappa\, d\hz\ .
\eeq
We now want to modify the proof of (\ref{me.16}) somewhat
to show that the above estimate holds for a.e. $t$ along a weak
$H^k$-flow. 

\bl\lb{me.l7} Let $N^3$ be complete, simply connected
manifold with sectional curvatures bounded above by $-\kappa$ for some
$\kappa>0$. Let $\Omega\subset N^3$ be a bounded, open set
with smooth boundary and $H|_{\p\Omega}>0$, which minimizes area from the outside in
$N^3$. Let $u$ be a weak $H^k$-flow generated by $\Omega$
and $k>2$. Then for all $t\in \tilde{B}$ we have
$$ \int_{\Gamma_t} \frac{1}{4}|\mathbf{H}|^2 -\kappa\, d\hz \geq 4\pi\
.$$
\el
\bpf The proof is  modification of Proposition \ref{me.p5}. We
first need the following claim.\\[+1.5ex]
\underline{Claim 1:} The sets $E_t$ minimize area from the outside in
$N^3$ for all $t \in (0,T)$.\\[+1.5ex]
Let $F \subset N^3$ be a comparison set for $E_{t}$, $t\in
(0,T)$, i.e. $E_t\subset F$. Now pick a time $\tau$ with
$t<\tau<T$. Then $E_\tau\Subset F$, thus we can pick a sequence of
sets $F_i$ with smooth boundary, such that $E_\tau \subset F_i$ for
all $i$ and $|\p F_i| \ra |\p^*F|$. By inequality (\ref{pr.9a}) we
have that $|\p^*(F_i\cap \Omega)| \leq |\p F_i|$. Since $\p \Omega$
and $\p F_i$ are smooth we can approximate $F_i\cap \Omega$ with sets
$F^i_j\Subset \Omega$ such that $E_\tau \subset F^i_j$ for all $j$ and
$|\p^*F^i_j|\ra |\p^*(F_i\cap\Omega)|$. Now the sets $F_j^i$ are valid
comparison sets for $E_\tau$ in $\Omega$, and thus taking the limits
$j\ra \infty$ and $i \ra \infty$ we see that
$$ |\p^*E_\tau|\leq |\p^* F|\ .$$
Now take a sequence $\tau_i \nearrow t$, then by non-fattening
$E_{\tau_i}\ra E_t$ and by outward minimizing $|\p^*E_{\tau_i}| \ra
|\p^*E_t|$, which proves the claim.\\[+1.5ex]
We will now show the corresponding result to Lemma \ref{me.l2}, with a
modification since the ambient space $N^3$ is not
flat. \\[1.5ex]
\underline{Claim 2:} Under the same conditions and with the same
notation as in Corollary \ref{me.c4} we have :
$$ \inl_{L_s(\Gamma_t)} f(y)\, d\hz(y) \leq  \inl_{\Gamma_t}
\tilde{f}(x)\big(\cosh(\gamma s) +
\frac{1}{2\gamma} \sinh(\gamma s) |\mathbf{H}(x)|\big)^2\,
d\hz(x)\ ,$$
where $\gamma =\sqrt{(C/2)},\ C:=\inf\{\, \bar{\text{R}}\text{ic}(X,X)\ |\ X \in T_pN^3,\
|X|=1,\ \text{dist}(p,\Gamma_t)\leq s \}$.\\[1.5ex]
The proof of this claim is nearly identical to the proof of Lemma \ref{me.l2}. The
only difference is that instead of (\ref{me.7}) we have
$$\frac{\p}{\p s} H = -|A|^2 - \bar{\text{R}}\text{ic}(\nu,\nu) \leq
-\frac{1}{2}H^2 + C\ , $$
which can be integrated to give
\beq\lb{me.35} H(s) \leq \frac{\sqrt{2C}\sinh\big(\sqrt{(C/2)}\cdot s\big) + 
\cosh\big(\sqrt{(C/2)}\cdot s\big) H(0)}{\cosh\big(\sqrt{(C/2)}\cdot
s\big)+ (\sqrt{2C})^{-1}\sinh\big(\sqrt{(C/2)}\cdot s\big) H(0)} \ .\eeq
Inserting this into (\ref{me.10}) we obtain
$$ d\mu(s) \leq \Big(\cosh\big(\sqrt{(C/2)}\cdot
s\big)+ (\sqrt{2C})^{-1}\sinh\big(\sqrt{(C/2)}\cdot s\big)
H(0)\Big)^2\, d\mu(0) \ .$$
Note that Lemma \ref{me.l3} remains true, so the rest of the proof follows again in the
same way.\\[+1.5ex]
We now have to rework Proposition \ref{me.p5}. Claim 1 in the proof
there remains true, but we have to replace the second claim by the
following.\\[+1.5ex]
\underline{Claim 3:} Under the same conditions and with the same
notation as in Claim 2 of the proof of Proposition \ref{me.p5}, we
have
$$ \frac{
\cosh\big(\sqrt{(C/2)}s\big) H_{\Gamma_t}(p)}{\cosh\big(\sqrt{(C/2)}
s\big)+ (\sqrt{2C})^{-1}\sinh\big(\sqrt{(C/2)}s\big)
H_{\Gamma_t}(p)}+ C s \geq H_{C_s}(\pi^{-1}_s(p))\\ .$$
As in the proof there let $\Sigma$ be a supporting hypersurface of
$C_s$ at $q$ such that (\ref{me.20}) holds. We may assume that
$\Sigma$ is a convex hypersurface, with principal curvatures locally
bounded by $1/s-\delta/4$. Since in $N^3$ there is no notion of
translating $\Sigma$ to touch $\Gamma_t$ in $p$ from the outside,
we directly evolve $\Sigma$ equidistantly towards $\Gamma_t$,
i.e. until $\Sigma_s$ touches $\Gamma_t$ in $p$. Since in an ambient Hadamard
manifold $N$, focal points develop later as in Euclidean
space the conditions above guarantee that $\Sigma_s$ remains
smooth. Thus
$$H_{\Gamma_t}(p)\geq H_{\Sigma_s}(p) $$
and by (\ref{me.35})
\beqs \begin{split}H_{C_s}-\delta\leq&\ H_{\Sigma}(q) \leq \frac{\sqrt{2C}\sinh\big(\sqrt{(C/2)}s\big) + 
\cosh\big(\sqrt{(C/2)} s\big) H_{\Sigma_s}(p)}{\cosh\big(\sqrt{(C/2)}
s\big)+ (\sqrt{2C})^{-1}\sinh\big(\sqrt{(C/2)} s\big)
H_{\Sigma_s}(p)}\\
\leq&\ \frac{\cosh\big(\sqrt{(C/2)} s\big) H_{\Sigma_s}(p)}{\cosh\big(\sqrt{(C/2)}
s\big)+ (\sqrt{2C})^{-1}\sinh\big(\sqrt{(C/2)} s\big)
H_{\Sigma_s}(p)} + C s\\
\leq&\ \frac{\cosh\big(\sqrt{(C/2)} s\big) H_{\Gamma_t}(p)}{\cosh\big(\sqrt{(C/2)}
s\big)+ (\sqrt{2C})^{-1}\sinh\big(\sqrt{(C/2)} s\big)
H_{\Gamma_t}(p)} + C s\ .
\end{split}\eeqs
which proves Claim 3 as $\delta \ra 0$.\\[+1.5ex]
To finally prove the lemma we combine (\ref{me.34}) and Claim 1
to estimate
\beqs\begin{split}
4\pi \leq&\ \inl_{C_s} G -\kappa\, d\hz = \!\!\!\inl_{C_s\cap
  L_s(\Gamma_t)}\!\!\!\!\!\!G\, d\hz+ \inl_{C_s\setminus
  L_s(\Gamma_t)}\!\!\!\!\!\!G\, d\hz - \inl_{C_s} \kappa\, d\hz\\ 
\leq&\ \!\!\!\inl_{C_s\cap L_s(\Gamma_t)}\!\!\!\!\!\!\frac{1}{4}H^2\,
d\hz +
\inl_{C_s\setminus
  L_s(\Gamma_t)}\!\!\!\!\!\!G\, d\hz - \inl_{\Gamma_t} \kappa\, d\hz
\end{split}\eeqs
By an argument of B. Kleiner, see the last part of the proof of
Proposition 8 in \cite{Kleiner92}, the second term in the second line goes
to zero as $s\ra 0$. So we can apply Claim 2 and 3 to estimate further
\begin{align*}
4 \pi\leq&\ \frac{1}{4}\!\!\!\!\inl_{\Gamma_t\cap C_0}\!\!\!\!\big(H_{C_s}(\pi^{-1}_s(x))\big)^2
\big(\cosh(\gamma s) +
\frac{1}{2\gamma} \sinh(\gamma s) |\mathbf{H}(x)|\big)^2\,
d\hz\\
&\ -\inl_{\Gamma_t} \kappa\, d\hz + o(s)\displaybreak[0] \\
\leq&\ \frac{1}{4}\!\inl_{\Gamma_t\cap C_0}\!\! (\cosh(\gamma
s))^2 |\mathbf{H}|^2\, d\hz - \inl_{\Gamma_t} \kappa\, d\hz + o(s) \\
\leq&\ \frac{1}{4}\inl_{\Gamma_t}\!\! 
|\mathbf{H}|^2\, d\hz - \inl_{\Gamma_t} \kappa\, d\hz + o(s) .
\end{align*}
Taking the limit $s\ra 0$ we obtain the estimate.\epf
 

\section{The monotonicity of the isoperimetric difference}\lb{monotonicity}

\noindent Let $\varphi \in C^1_\text{c}(\R),\ \varphi\geq 0$ such that $\int_\R
\varphi \, dx \geq 1$. Choose a function $\phi:= \varphi(x_{n+2})
\in C^1(\Omega\times \R)$ to  define the approximative area and volume by
\beq\lb{mo.00} A^i_t:= \int_{N^i_t} \varphi\, d\hne\ ,\ \ \ \ V^i_t:= \int_{E^i_t}
\varphi \, d\hnz 
\eeq
and the approximate isoperimetric difference 
$$I^i_t:= \big(A^i_t\big)^\frac{n+1}{n}- c_{n+1}V^i_t\ $$
where $c_{n+1}$ is defined as in the introduction. Let $t_1,t_2 \in
(0,T)$. Assuming that $i$ is big enough, the boundary (in
$N\times \R$) of the graphs $N^i_t$ does not intersect
supp$(\varphi)$ for $t \in (t_1,t_2)$. By the Coarea formula and the
evolution equations we have
\beqs
V^i_{t_2}-V^{i}_{t_1} =  - \int_{t_1}^{t_2}\!\int \varphi H_i^k\,
d\mu^i_t\, dt
\eeqs
and
\beqs\begin{split}
\big(A^i_{t_2}\big)^\frac{n+1}{n}-\big(A^i_{t_1}\big)^\frac{n+1}{n}
=&\ 
- \frac{n+1}{n}\int_{t_1}^{t_2}\big(A^i_t\big)^\frac{1}{n}\int
\varphi H_i^{k+1}\, d\mu^i_t\, dt\\
&\ \ - \frac{n+1}{n}\int_{t_1}^{t_2}\big(A^i_t\big)^\frac{1}{n}\int
\langle \nabla\varphi, \nu\rangle  H_i^{k}\, d\mu^i_t\, dt\ .
\end{split}\eeqs
By H\"older's inequality we can estimate
\begin{align}\lb{mo.0}
\frac{I^i_{t_2}-I^i_{t_1}}{c_{n+1}} \leq \ \int_{t_1}^{t_2}&\!\!\Bigg(\!\bigg(\int \varphi
H_i^{k+1}\,
d\mu^i_t\bigg)^\frac{k}{k+1}\!\!\big(A^i_t\big)^\frac{1}{k+1} \\ &\ - 
\frac{n+1}{n c_{n+1}}\big(A^i_t\big)^\frac{1}{n}\int
\varphi H_i^{k+1}\, d\mu^i_t\Bigg) dt\displaybreak[1]\\
&\ -\frac{n+1}{nc_{n+1}}\int_{t_1}^{t_2}\big(A^i_t\big)^\frac{1}{n}\int
\langle \nabla\varphi, \nu\rangle  H_i^{k}\, d\mu^i_t\, dt\ .
\end{align}

\bl\lb{mo.l1} Let $k\geq n-1$. Assume that for a.e. $t \in [0,T]$ the estimate 
\beq\lb{mo.1} 
\int_{\Gamma_t}|\mathbf{H}|^n\, d\mu \geq \Big(\frac{n}{n+1}\,
c_n\Big)^n
\eeq holds. Then for any $t_1,t_2 \in (0,T)$, $t_1<t_2$,  
\beqs
\limsup_{i\ra \infty} \int_{t_1}^{t_2}\!L^i_t\, dt\leq 0\ ,
\eeqs
where
\beqs
L^i_t:=\bigg(\int \varphi
H_i^{k+1}\,
d\mu^i_t\bigg)^\frac{k}{k+1}\!\!\big(A^i_t\big)^\frac{1}{k+1} -
\frac{n+1}{n c_{n+1}}\big(A^i_t\big)^\frac{1}{n}\int
\varphi H_i^{k+1}\, d\mu^i_t\ .
\eeqs
\el
\bpf Since $t_2<T$ and $\phi \not \equiv 0$ there is $\delta >0$
such that $A^i_t \geq \delta$ for all $t \in [t_1,t_2]$ and all
$i$. Since $k/(k+1)<1$, there is
a $C_1\geq 0$ such that
\beqs
L^i_t \leq C_1
\eeqs
for all $t\in [t_1,t_2]$ and all $i$. To prove the lemma it thus
suffices by Fatou's lemma to show that $\limsup_{i\ra\infty} L^i_t \leq 0$ for all
$t\in \tilde{B}\cap [t_1,t_2]$ such that (\ref{mo.1}) holds. Fix such a $t$. Arguing as above we
see that there is a $C_2 >0$, independent of $i$, such that
$$ L^i_t \leq 0 \ \ \text{if}\ \ \int \varphi H_i^{k+1}\, d\mu^i_t
\geq C_2\ .$$
Thus we can assume that 
$$ \int \varphi H_i^{k+1}\, d\mu^i_t
\leq C_2\ \ \ \text{for all}\ i \ .$$
We write $L^i_t$ in the form $L^i_t= a_i \cdot b_i$, where
\beqs a_i:=\bigg(\big(A^i_t\big)^\frac{1}{k}\int \varphi
H_i^{k+1}\,
d\mu^i_t\bigg)^\frac{k}{k+1}
\eeqs
and
\beqs
b_i:=1 -
\frac{n+1}{n c_{n+1}}\big(A^i_t\big)^{\frac{1}{n}-\frac{1}{k+1}}\bigg(\int
\varphi H_i^{k+1}\, d\mu^i_t\bigg)^\frac{1}{k+1}
\eeqs
with the bounds $0\leq a_i \leq C_2^\prime\ ,\ |b_i|\leq C_3$\ . By
the lower semicontinuity of the $L^{k+1}$-norm of $H^i$ and
(\ref{mo.1}) we see that
$$\limsup_{i\ra \infty} b_i \leq 1 -
\frac{n+1}{n c_{n+1}}\big(\mu_t(\varphi)\big)^{\frac{1}{n}-\frac{1}{k+1}}\bigg(\int
\varphi |\mathbf{H}|^{k+1}\, d\mu_t\bigg)^\frac{1}{k+1} \leq 0\ , $$
since $\mu_t=\hne\res(\Gamma_t\times\R)$ and
$\varphi(x)=\phi(x_{n+2})$ with $\int_\R \phi\, dx \geq 1$\ . Since
the $a_i$ are nonnegative and uniformly bounded this implies also
$\limsup_{i \ra \infty} L^i_t \leq 0$.
\epf

\bp \lb{mo.p1} Let $k\geq n-1$ and $u$ be a weak $H^k$-flow generated by
$\Omega$. Assume that $(\ref{mo.1})$ holds for a.e. $t \in
[0,T]$. Then the isoperimetric difference
$$ I_t:= \big(\hn(\p^*E_t)\big)^\frac{n+1}{n} - c_{n+1}\hne(E_t) $$
is monotonically decreasing along $u$ for $t \in [0,T)$.
\ep
\bpf Let $l>1$. Pick a $\phi\in C^1_c(\R),\ 0\leq\phi\leq 1/(2l)$ with $\phi =
1/(2l)$ on $[-l,l]$, $\phi = 0$ on $\R\setminus [-l-1,l+1]$, and
$|D\phi|\leq 1/l$. By the mass bound we have
\beqs
\Big|\int_{t_1}^{t_2} (A^i_t)^\frac{1}{n}\int \langle \nabla\phi,\nu
\rangle H_i^k\, d\mu_t^i\, dt \Big| \leq \frac{C_3}{l}\ .
\eeqs
For $t \in \tilde{B}$ we have
$$ I^i_t \ra \bigg(\int_{\Gamma_t\times \R} \varphi\,
d\hne\bigg)^\frac{n+1}{n}- c_{n+1} \int_{E^\prime_t} \varphi\, d\hnz\
:= I_t^\prime.$$
By the above estimate, equation (\ref{mo.0}) and Lemma
\ref{mo.l1} we can deduce that in the limit
$$ I^\prime_{t_2} \leq I^\prime_{t_1} + \frac{C_3}{l} $$
for any $t_1,t_2 \in \tilde{B},\ t_1< t_2$. Since $I^\prime_t \ra I_t$
as $l\ra \infty$ this implies that
$$ I_{t_2}\leq I_{t_1}\ . $$
Since $u$ is non-fattening, and all sets $E_t$ minimize area from the
outside in $\Omega$, we can approximate times in
$[0,T]\setminus \tilde{B}$ with times in $\tilde{B}$ to see that this
monotonicity holds for all $t\in [0,T)$. 
\epf

\noindent {\it Proof of Theorem \ref{mo.t1}:} By Lemma \ref{me.l6} and Proposition \ref{me.p5} together with
Proposition \ref{mo.p1} the isoperimetric difference is monotonically
decreasing in all of the above cases. Since $u$ is non-fattening we have $\lim_{t\ra T}I_t \geq 0$.
\hfill $\qedbox$\\[1.5ex]
We can now use the above theorem to a give a proof of the
isoperimetric inequality.\\[+1.5ex]
 \noindent {\it Proof of Corollary \ref{mo.c1}:} We can first assume that $U$ is connected. Let $\Omega$ be the
outer minimizing hull of $U$, see the proof of Theorem \ref{eu.t1}
for more details on minimizing hulls. Since $U\subset \Omega$ and
$|\p^*\Omega|\leq |\p U|$ we can replace $U$ by $\Omega$. We first
treat the case $n \leq 6$. Here we can apply a result of
Sternberg, Williams and Ziemer \cite{SternbergWilliamsZiemer91}, which
shows that $\p\Omega$ is a $C^{1,1}$-hypersurface. Thus $\p\Omega$
carries a weak mean curvature in $L^\infty$, which is non-negative,
since $\Omega$ minimizes area from the outside. We can thus apply a
result of G. Huisken and T. Ilmanen, Lemma 2.5 in
\cite{HuiskenIlmanen02}, which states that starting from such a
hypersurface there exists a smooth solution to mean curvature flow
$(M_t)_{0< t< \eps}$, where all the $M_t$ have strictly positive mean
curvature for $t>0$. Furthermore, the initial datum $M_0=\p\Omega$ is
attained in $C^{1,\alpha}$ and since the mean curvature is positive,
the hypersurfaces $M_t$ foliate a neighborhood of $\p\Omega$ in
$\Omega$. The lemma in \cite{HuiskenIlmanen02} is stated only if the
ambient space is flat, but by a closer examination one sees that with
some minor modifications the same proof also works in $N^3$
as above. Then for $0<t<\eps$ let $\Omega_t:= \Omega 
\setminus \big(\cup_{0<\tau\leq t}M_t\big)$, thus $\p\Omega_t =M_t$. Pick any
$t\in (0,\eps)$.
Since $\p\Omega_t$ has strictly positive
mean curvature, there exists a weak $H^k$-flow, generated by
$\Omega_t$. Taking $k>n$ if $n\geq 4$, or $k\geq 1$ in the case $n=3$,
we can apply Theorem \ref{mo.t1} to see that $\Omega_t$ satisfies
(\ref{mo.2}). We then take the limit $t\searrow 0$.\\[1ex]
For $n=7$ the outer minimizing hull can have isolated singularites
in the part where it does not touch the obstacle, i.e. in the part
where its boundary constitutes a minimal surface. To treat this case
let 
$$U_\tau:=\{x\in U\ |\ \text{dist}(x,\p U)> \tau\},$$ which has a
smooth boundary for small enough $\tau>0$. Let $E$ be the outer
minimzing hull of $U$ and $E_\tau$ be the outer minimizing hulls of
$U_\tau$. We have $E_s\subset E_\tau\subset E$ for $0<\tau<s$, which
implies by Lemma \ref{pr.l5} that $E_\tau \ra E$ in $L^1$ as $\tau \searrow 0$ as well as 
$\p^*E_\tau \ra \p^*E$ in the sense of radon measures. We now apply a
strong maximum principle of Ziemer and Zumbrun \cite{ZiemerZumbrun99}
(which actually is an application of the strong maximum principle of Moschen-Simon) to show
that $E_\tau \Subset E$ for $\tau>0$. This maximum principle states
that if two sets $F,G$ are outer minimizing and minimizing, respectively, with respect
to an open set $V$ and $G\cap V \subset F\cap V$, then either $\p F =
\p G$ in $V$ or $\p F\cap\p G =\emptyset$ relative to $V$, provided
$\p F\cap V$ and $\p G \cap V$ are connected. Note that
$E_\tau$ minimizes area in $N\setminus U_\tau$ and that $\p E$ is
connected in $N\setminus U_\tau$ and does not touch
$U_\tau$. Thus applying the maximum principle to every connected
component of $\p E_\tau \cap (N\setminus U_\tau)$ we see that $E_\tau \Subset
E$ for all $\tau > 0$.\\
These properties enable us to argue as Hardt-Simon in
\cite{HardtSimon85}, Theorem 5.6, to deduce that $\p E_\tau$ is
$C^{1,1}$ for all $\tau$ small enough. Note that $\p E$ is
$C^{1,1}$ in a neighbourhood of $\p U$, and thus also $\p E_\tau$ for
small enough $\tau$. Thus we can replace $U$ by $U_\tau$ and argue as
for $n\leq 6$, finally taking the limit $\tau \searrow 0$.  \hfill $\qedbox$\\[1ex]

In the case that $N^3$ has sectional curvatures bounded
above by $-\kappa$, $\kappa>0$, we aim to apply estimate (\ref{me.34})
to show how one can use a
weak $H^k$-flow to prove that the
isoperimetric profile of $N$ always lies above the
isoperimetric profile of the model space with constant curvature
$-\kappa$.\\
Let $(M_t)_{0\leq t <T}$ be a smooth $H^k$-flow in $N^3$ of
hypersurfaces with positive mean curvature. Let us assume all the
$M_t$ are diffeomorphic to a sphere. We then can apply (\ref{me.34})
to estimate 
\beqs \begin{split}
-\frac{d}{dt}V =&\ \inl_{M_t} H^k\, d\hz \leq \bigg( \inl_{M_t}
H^{k+1}\, d\hz\bigg)^\frac{k}{k+1} A^\frac{1}{k+1}\\
&\ \cdot (16\pi+4\kappa A)^{-\frac{1}{2}}\bigg( \inl_{M_t}
H^{k+1}\, d\hz\bigg)^\frac{1}{k+1} A^{\frac{1}{2}-\frac{1}{k+1}}\\
=&\ (16\pi+4\kappa A)^{-\frac{1}{2}}A^\frac{1}{2}\inl_{M_t}
H^{k+1}\, d\hz = -\frac{d}{dt}f_\kappa(A)\ ,
\end{split} \eeqs
where $f_\kappa$ is defined by (\ref{intro.3}). 
Thus $f_\kappa(A)-V$ is monotonically decreasing under the flow. Consider the case
that $N_\kappa^3$ is the model space of constant curvature
$-\kappa$ and let $M_t$ be the $H^k$-flow of geodesic spheres
contracting to a point. Then (\ref{me.34}) holds with equality for all
$M_t$ and also the above calculation is an equality. Using that in the
model space geodesic balls optimize the isoperimetric ratio, we have 
 $$\mathcal{H}^3(U) \leq f_\kappa(\hz(\p U))\ ,$$ 
for all open and bounded $U\subset N^3_\kappa$, with
equality on geodesic balls. Arguing as in the beginning of this
section, and using Lemma \ref{me.l7} we arrive at the following
proposition.

\bp\lb{mo.p2} Let $k>2$ and $u$ be a weak $H^k$-flow generated by
$\Omega$, where $\Omega\subset N^3$ is open and
bounded, and $N^3$ has sectional curvatures bounded above by
$-\kappa$. Furthermore assume that $\p\Omega$ is smooth with strictly
positive mean curvature and that $\Omega$ minimizes area from the
outside in $N^3$. Then
$$I_t^\kappa:=f_\kappa(\hz(\p^*E_t)) - \mathcal{H}^3(E_t)$$
is a nonnegative, monotonically decreasing function along $u$ for $t\in [0,T)$.
\ep
\bpf
We define $A^i_t$ and $V^i_t$ as in (\ref{mo.00}). Then take
$$ I_t^{\kappa,i}:=f_\kappa(A^i_t)-V^i_t\ .$$ 
Arguing as before we can estimate
\beqs\begin{split}
I^{\kappa,i}_{t_2}-I^{\kappa,i}_{t_1}\leq& \ \int_{t_1}^{t_2}\!\!\Bigg(\!\bigg(\int \varphi
H_i^{k+1}\,
d\mu^i_t\bigg)^\frac{k}{k+1}\!\!\big(A^i_t\big)^\frac{1}{k+1} \\ -&\ 
\frac{\big(A^i_t\big)^\frac{1}{2}}{\big(16\pi + 4\kappa A^i_t\big)^\frac{1}{2}}\int
\varphi H_i^{k+1}\, d\mu^i_t\Bigg) dt\\
&\ -\int_{t_1}^{t_2}\frac{\big(A^i_t\big)^\frac{1}{2}}{\big(16\pi + 4\kappa A^i_t\big)^\frac{1}{2}}\int
\langle \nabla\varphi, \nu\rangle  H_i^{k}\, d\mu^i_t\, dt\ .
\end{split} \eeqs
Analogous to the proof of Lemma \ref{mo.l1} we can use Lemma
\ref{me.l7} to show that for $t_1,t_2 \in (0,T),\ t_1<t_2$, we have
\beqs\begin{split}\limsup_{i\ra \infty}\inl_{t_1}^{t_2} \Bigg(\!\bigg(\int \varphi
H_i^{k+1}\,
&\ d\mu^i_t\bigg)^\frac{k}{k+1}\!\!\big(A^i_t\big)^\frac{1}{k+1}\\
 &\ -\frac{\big(A^i_t\big)^\frac{1}{2}}{\big(16\pi + 4\kappa A^i_t\big)^\frac{1}{2}}\int
\varphi H_i^{k+1}\, d\mu^i_t\Bigg)\, dt \leq 0\ . 
\end{split}\eeqs
The rest follows as in Proposition \ref{mo.p1} and Theorem \ref{mo.t1}.
\epf

This enables us to also give a new proof of the stronger result 
in the case that the sectional curvatures of
$N^3$ are bounded above by $-\kappa <0$.\\[2ex]
\noindent {\it Proof of Theorem \ref{mo.t2}:} We can use Proposition \ref{mo.p2} to give an analogous proof as
in Corollary \ref{mo.c1}. We use the same terminology as in the proof
there. The only missing bit to apply
Proposition \ref{mo.p2} is to show that the sets
$\Omega_t$ are minimizing area from the outside in $N^3$. Fix
a $t>0$. Take $F$ to be a comparison set for
$\Omega_t$, i.e. $\Omega_t \subset F$. Again by (\ref{pr.9a}), we have
that $|\p^*(F\cap \Omega)| \leq |\p^*F|$ since $\Omega$ minimizes area
from the outside in $N^3$. Note
that the surfaces $M_\tau =\p \Omega_\tau,\ \tau\in (0,t)$ smoothly foliate
$\Omega\setminus \Omega_t$ and $M_t \ra \p\Omega$ in $C^1$ as $t\ra
0$. Since all $M_t$ have nonnegative mean curvature, we can use the outer
unit normal vectorfield $\nu$ to these hypersurfaces as a calibration on
$\Omega\setminus \Omega_t$ which satisfies
$\text{div}_{N^3}(\nu)\geq 0$. Using this calibration we see
that $|\p\Omega_t| \leq |\p^*(F\cap \Omega)|$, i.e. $\Omega_t$
minimizes area from the outside in $N^3$. \hfill $\qedbox$

\bibliographystyle{plain}

\providecommand{\bysame}{\leavevmode\hbox to3em{\hrulefill}\thinspace}
\providecommand{\MR}{\relax\ifhmode\unskip\space\fi MR }
\providecommand{\MRhref}[2]{%
  \href{http://www.ams.org/mathscinet-getitem?mr=#1}{#2}
}
\providecommand{\href}[2]{#2}

\end{document}